\documentclass[11pt]{amsart}
\usepackage[T1]{fontenc}
\usepackage[utf8]{inputenc}
\usepackage{lmodern}
\usepackage[a4paper,margin=29mm]{geometry}
\usepackage{amsmath,amssymb,amsthm,mathtools,bm}
\usepackage{enumitem}
\usepackage{tikz}
\usetikzlibrary{calc,arrows.meta}
\usepackage[unicode,hidelinks]{hyperref}
\newtheorem{theorem}{Theorem}[section]
\newtheorem{proposition}[theorem]{Proposition}
\newtheorem{lemma}[theorem]{Lemma}
\newtheorem{corollary}[theorem]{Corollary}
\theoremstyle{definition}

\theoremstyle{remark}

\newcommand{\one}{\mathbf 1}

\DeclareMathOperator{\Per}{Per}

\numberwithin{equation}{section}
\allowdisplaybreaks[2]
\title[Triple bubbles in the hyperbolic plane]{Triple bubbles in the hyperbolic plane}
\author{Zihao Wang}
\address{School of Mathematical Sciences, Fudan University, Shanghai 200433, China}
\email{wangzh25@m.fudan.edu.cn}
\date{1 November 2026}
\begin{document}
\raggedbottom
\begin{abstract}
We prove that the standard triple bubble uniquely minimizes perimeter
among clusters in the hyperbolic plane with three prescribed positive
finite areas.  The standard family admits global pressure coordinates,
and its volume differential is the negative of a positive definite
grounded weighted graph Laplacian.  A normalized perimeter excess attains a comparison minimum
under a perimeter bound and gives a strict comparison between actual
and model pressures.  Component relabeling reduces the possible
topologies to an all-boundary network or a 4-daisy.  Discrete
comparisons exclude the 4-daisy; pressure-preserving face deletion
reduces the remaining network to the standard tetrahedral configuration.
\end{abstract}
\maketitle

\section{Introduction}

The multi-bubble problem asks for the least perimeter needed to
enclose and separate several prescribed volumes.  Even on a surface,
regularity alone leaves a substantial global problem: constant-
curvature edges and $120$-degree junctions admit many topologies,
and a single cell label may occupy several connected components.
In the hyperbolic plane there is a further obstruction to the usual
scaling argument, since dilation changes the curvature of the
ambient metric.

Throughout the paper the hyperbolic plane has sectional curvature
$-1$. We write $|E|$ for hyperbolic area and $\mathcal H^d$ for
$d$-dimensional Hausdorff measure in the ambient metric. For an open
set $U\subset\mathbb H^2$, the perimeter of a measurable set $E$ in
$U$ is
\[
 \Per(E;U)=\sup\left\{\int_E\operatorname{div}_H X\,dV_H:
 X\in C_c^1(U,T\mathbb H^2),\ |X|\le1\right\}.
\]
The set has locally finite perimeter if this quantity is finite for
every relatively compact $U$; we put $\Per(E)=\Per(E;\mathbb H^2)$.
The reduced boundary $\partial^*E$ and its measure-theoretic outward
normal are understood in the sense of \cite[Chapters 12 and 15]{Maggi};
in particular, the perimeter measure is
$\mathcal H^1\!\llcorner\partial^*E$.
A three-cluster is a measurable partition
$\mathcal C=(E_0,E_1,E_2,E_3)$, up to null sets, into cells of
locally finite perimeter, where
\[
 |E_i|=v_i\in(0,\infty)\quad(1\le i\le3),\qquad
 E_0=\mathbb H^2\setminus(E_1\cup E_2\cup E_3).
\]
Its perimeter is
\begin{equation}\label{eq:cluster-perimeter}
 P(\mathcal C)=\frac12\sum_{i=0}^3\Per(E_i)
 =\sum_{0\le i<j\le3}
       \mathcal H^1(\partial^*E_i\cap\partial^*E_j).
\end{equation}
Thus every interface is counted once. We minimize over all such
partitions, with no connectedness or topological constraint on the
cells. A cluster has finite perimeter when $P(\mathcal C)<\infty$;
we use BV for functions of bounded variation.
The infinite-area cell $E_0$ is the exterior; $E_1,E_2,E_3$
are the finite cells.

A standard hyperbolic triple bubble is defined as follows. Identify
$\mathbb H^2$ conformally with
$\mathbb S^2_+=\{x\in\mathbb S^2:x_3>0\}$, equipped with
$x_3^{-2}g_{\mathbb S^2}$. Choose unit vectors
$a_0,a_1,a_2,a_3\in\mathbb R^3$ satisfying
\[
 \sum_{i=0}^3a_i=0,\qquad a_i\cdot a_j=-\frac13\quad(i\ne j),
 \qquad U_i=\{x\in\mathbb S^2:a_i\cdot x<a_j\cdot x
                                  \text{ for all }j\ne i\}.
\]
For a M\"obius transformation $\Phi$ of $\mathbb S^2$ with
$\overline{\Phi(U_i)}\subset\mathbb S^2_+$ for $1\le i\le3$,
the cells $\Phi(U_i)\cap\mathbb S^2_+$ form a standard triple bubble.
Here a M\"obius transformation is a conformal automorphism of the
round sphere, allowing orientation reversal. The inequalities defining
$U_i$ fix both the interfaces and the sides assigned to their labels.
Section~\ref{sec:standard} constructs this family in pressure
coordinates and proves that its labeled areas determine it uniquely
up to hyperbolic isometry. Let $T_k$ denote the standard $k$-bubble
perimeter profile, and write $T=T_3$.
The minimizing profile is
\[
 I_3(v)=\inf\{P(\mathcal C):|E_i|=v_i,\ 1\le i\le3\}.
\]

\begin{theorem}\label{thm:main}
For every $v\in(0,\infty)^3$ and every finite-perimeter
three-cluster of finite-cell areas $v$,
\[
 P(\mathcal C)\ge T(v).
\]
Equality holds precisely for the standard triple bubble, up to
hyperbolic isometry and null-set modifications, with the cell labels
retained.
\end{theorem}

Foisy, Alfaro, Brock, Hodges, and Zimba proved the Euclidean planar
double-bubble theorem \cite{FoisyEtAl}. Hutchings, Morgan, Ritor\'e, and
Ros proved its three-dimensional counterpart \cite{HMRR}, and
Reichardt extended the theorem to every Euclidean dimension
$n\ge3$ \cite{Reichardt}. Cotton--Freeman proved the double-bubble theorem
in the hyperbolic plane \cite[Theorem 2.7]{CottonFreeman}.
For three prescribed areas, Wichiramala proved the Euclidean planar
triple-bubble theorem \cite{Wichiramala}. Lawlor gave a common
unification argument for the Euclidean and spherical planar
triple-bubble problems \cite{Lawlor}. Milman--Neeman proved the
Euclidean and spherical multi-bubble conjectures for
$2\le k\le\min\{4,n\}$ finite bubbles
\cite[Theorem 1.13]{MilmanNeemanStructure}. Milman's survey discusses the
hyperbolic planar triple-bubble problem as a remaining case
\cite{MilmanSurvey}.

\subsection*{Outline of the proof}

We adapt Lawlor's discrete unification method to the hyperbolic
plane. For $A(v)=v_1+v_2+v_3$, define
\begin{equation}\label{eq:ratio}
 S(v)=T(v)-A(v),\qquad
 J(\mathcal C)=P(\mathcal C)-A(v),\qquad
 Q(\mathcal C)=\frac{J(\mathcal C)}{S(v)}.
\end{equation}
Subtracting total area removes the leading linear growth of
hyperbolic perimeter. Both excesses are positive, and a strict
counterexample satisfies $0<Q<1$.

The standard family has global pressure coordinates
$\kappa\in(1,\infty)^3$. Its volume differential is the negative
of a positive definite grounded weighted graph Laplacian, and
\[
 D_\kappa v=-L_w,\qquad
 dA_{\mathrm{std}}=-\sum_iw_{i0}\,d\kappa_i,
 \qquad w_{i0}>0.
\]
In particular, decreasing all pressures strictly increases total
area. If a strict counterexample exists, uniform small-area and
vanishing-cell estimates yield a comparison minimum below one under
a perimeter bound.
First variation at this minimum gives $0<p_i<\kappa_i$ for the
cluster pressures $p_i$.

Component relabeling and the 4-daisy comparison leave a network
in which every finite face meets the exterior. Successive local
face deletions preserve the pressures and reduce this network to a
standard triple bubble; each deletion strictly decreases total
finite area. If no deletion occurs, the original cluster is standard
and has $Q=1$. Otherwise, the standard endpoint implies $p_i>1$. The resulting
inequalities $1<p_i<\kappa_i$ force the opposite area comparison
by the pressure differential above, giving a contradiction. Once
$Q\ge1$ is established, first variation at equality gives $p=\kappa$;
uniqueness of pressure realization then rules out any deletion.

\section{Minimizers and planar networks}\label{sec:preliminaries}

We use the open representatives of the cells of a minimizing cluster.  An
interface is oriented by the unit normal pointing from its first cell into
its second cell.  Writing $\nabla^H$ for the hyperbolic Levi--Civita connection, our
curvature convention is
\[
 \nabla^H_TN=kT.
\]
Thus the boundary of a hyperbolic disk of radius $R$, with its outward
normal, has curvature $\coth R$.  With exterior pressure $p_0=0$, the first
variation formulas are
\begin{equation}\label{eq:first-variation}
 \delta P=\sum_{i=1}^3p_i\,\delta v_i,
 \qquad k_{ij}=p_i-p_j.
\end{equation}

\begin{proposition}[Existence, boundedness, and regularity]\label{prop:inputs}
For every $v\in(0,\infty)^3$ the infimum $I_3(v)$ is attained.
Every minimizing cluster has bounded finite cells and connected boundary.
After modification on a null set, its boundary is a finite embedded
network of smooth arcs of constant geodesic curvature.  At each vertex
exactly three arcs meet at angles $2\pi/3$, and the three intervening
sectors have distinct cell labels.  The profile $I_3$ is continuous on
$(0,\infty)^3$.
\end{proposition}

\begin{proof}
The existence, boundedness, and local H\"older continuity results of
Resende de Oliveira \cite{Resende} apply to $\mathbb H^2$, since
\[
 \operatorname{Ric}_g=-g,\qquad
 \operatorname{inj}(\mathbb H^2)=\infty,\qquad
 |B_g(x,1)|=2\pi(\cosh 1-1)>0.
\]
Moreover, homogeneity identifies every pointed manifold
$(\mathbb H^2,g,x)$ isometrically with a fixed pointed copy of
$\mathbb H^2$.  Thus the pointed limits required in the existence
theorem are again the space form of curvature $-1$.  In particular,
for every $v\in(0,\infty)^3$ there are $C_v,\delta_v>0$ such that
\[
 |I_3(w)-I_3(v)|\le C_v|w-v|^{1/2}
 \qquad (w\in(0,\infty)^3,\ |w-v|<\delta_v).
\]

Let $\mathcal C=(E_0,E_1,E_2,E_3)$ be a minimizer. To apply
Morgan's surface regularity argument
\cite[Lemma 2.2, Theorem 2.3, Steps 2--3, and Section 2.4]{Morgan},
choose bounded representatives and concentric hyperbolic balls
$U\Subset U'$ with $\overline E_1\cup\overline E_2\cup\overline E_3
\subset U$.  Extend the metric on a
neighborhood of $\overline{U'}$ to a smooth metric on a compact
surface diffeomorphic to the sphere.  Extend the partition by putting
the complement of $U'$ in $E_0$.  All four cells have positive area,
and their boundary currents are rectifiable.  The reduced-boundary
construction in Morgan's Lemma 2.2 supplies paired area-adjustment supports within $U'$,
before any regularity of the interfaces is assumed.  Its conclusion
is the following: for some $C,\varepsilon>0$, every
$a=(a_1,a_2,a_3)$ with $|a|<\varepsilon$ can be realized by a smooth
diffeomorphism $\Phi_a$, with
\begin{equation}\label{eq:rough-area-adjustment}
 |\Phi_a(E_i)|-|E_i|=a_i\quad(1\le i\le3),\qquad
 P(\Phi_a\mathcal C)\le P(\mathcal C)+C|a|.
\end{equation}
The supports come in three pairs whose doubled balls are mutually
disjoint; either member of each pair may be used.  Hence there is
$r_0>0$ such that every ball $B_g(x,r)$ with $x\in U$ and $r<r_0$
lies in $U'$ and misses one complete choice of adjustment supports.

For a cluster $\mathcal F=(F_0,F_1,F_2,F_3)$ with
$F_i\mathbin\triangle E_i\Subset B_g(x,r)$, put
\[
 a_i=|E_i|-|F_i|
     =\int_{B_g(x,r)}(\one_{E_i}-\one_{F_i})\,dA_g.
\]
Choose the adjustment supports outside this ball.  On those supports
$\mathcal F=\mathcal C$, so the same diffeomorphism in
\eqref{eq:rough-area-adjustment} corrects $\mathcal F$ and has the
same perimeter cost.  Minimality, followed by $|a_i|\le|B_g(x,r)|$,
gives, after decreasing $r_0$ if necessary,
\begin{equation}\label{eq:local-perimeter-comparison}
 \begin{split}
 P(\mathcal C;B_g(x,r))
 &\le P(\mathcal F;B_g(x,r))+C|a|\\
 &\le P(\mathcal F;B_g(x,r))+\Lambda r^2.
 \end{split}
\end{equation}
Here $P(\mathcal C;B)$ denotes the sum of interface lengths in $B$;
$\Lambda$ is independent of $x\in U$ and $r<r_0$.

The remaining hypotheses of Morgan's Steps 2--3 are local on $U'$.
In fixed plane coordinates choose $c\ge1$ such that
$c^{-2}g_{\mathrm{eucl}}\le g\le c^2g_{\mathrm{eucl}}$ on $U'$.
For a sufficiently short rectifiable Jordan cycle $\gamma\subset U$,
the planar isoperimetric estimate used in Step 2 gives
\[
 |\operatorname{int}\gamma|_g
 \le c^2|\operatorname{int}\gamma|_{\mathrm{eucl}}
 \le\frac{c^2}{4\pi}\operatorname{length}_{\mathrm{eucl}}(\gamma)^2
 \le\frac{c^4}{4\pi}s^2,\qquad
 s=\operatorname{length}_g(\gamma).
\]
The cycle and its interior lie in $B_g(x,s)$ for $x\in\gamma$:
the shorter arc gives $\gamma\subset\overline{B_g(x,s/2)}$, and
hyperbolic balls are convex disks.  For $s<r_0$, the filling
operation in Step 2 saves at least $s/3$ of interface length;
\eqref{eq:rough-area-adjustment} corrects its area errors outside
this ball.  Consequently minimality gives
\[
 0\le-\frac{s}{3}+C_1s^2.
\]
This excludes arbitrarily short boundary cycles.  The finite-cycle
and separation-and-rejoining arguments in those steps therefore
apply on $U'$, as do the first-variation and tangent-network
arguments in Section 2.4.  They give a neighborhood of every
boundary point containing either one embedded smooth arc of
constant geodesic curvature or three such arcs meeting at $120$
degrees.

The boundary is compact.  Its arc and triple-junction neighborhoods
have a finite subcover, so it has finitely many vertices and edges.
The three sectors at a vertex have distinct labels, since each of
the three incident arcs separates two different cells.  Finally,
$\mathbb H^2$ is smooth, connected, simply connected, homogeneous,
and one-ended.  Milman--Neeman's connectedness theorem
\cite[Theorem 1.1]{Connected} therefore makes this boundary connected.
\end{proof}

A regular network is a finite embedded graph whose edges are smooth
arcs with disjoint relative interiors and whose vertices have three
incident arcs, three sectors of angle $2\pi/3$, and three distinct
cell labels. A regular cluster is a partition with such a boundary.
It is an equilibrium network if there is one global pressure per
cell label, with $p_0=0$, such that every oriented $i/j$ edge has constant
curvature $p_i-p_j$. Thus equilibrium is a stationarity condition;
fixed-area minimality is an additional property.
We call a cluster bounded when the union of its finite cells is
bounded, and always use compactness of a network to mean compactness
of its boundary graph.

A face of the network is a connected component of the complement of the
boundary graph, with its inherited cell label. A disk face is a
Jordan domain. Its number of sides counts the edges of its boundary
cycle, retaining parallel edges. A digon has two sides, and a
triangular face has three. A finite-cell component $K$ is a boundary
component if
$\mathcal H^1(\partial^*K\cap\partial^*E_0)>0$; otherwise it is an
interior component. An all-boundary network is one whose finite
faces are all boundary components.

\begin{lemma}[Faces of the minimizing network]\label{lem:disk-faces}
The boundary graph of a minimizing three-bubble has no loops,
bridges, or cut vertices.  On the one-point compactification its
faces have simple boundary cycles, including two-edge cycles when
parallel edges occur.  Every bounded face is a Jordan disk.  Each
face is one connected component of a cell.
\end{lemma}

\begin{proof}
The preceding proposition supplies a compact connected finite
trivalent graph.  A vertex-free circular component would be the
whole connected boundary and would leave only one bounded face;
three positive finite cells exclude this possibility.
A bridge in an embedded plane graph has the same complementary
face on its two sides, whereas an interface separates distinct
labels.  Thus there are no bridges.  A loop at a trivalent vertex
uses two half-edges, making its remaining edge a bridge.  If a
trivalent graph has a cut vertex, one component after its removal
is attached by just one of the three half-edges, again producing
a bridge.  The graph is therefore two-connected.  Its face
boundary walks are simple cycles: a repeated vertex in such a
walk separates two nonempty groups of incident edges and would
be a cut vertex.  The Jordan curve theorem identifies the bounded
complementary faces with disks.  Distinct faces bearing the same
label are distinct connected components of the corresponding open cell.
If $V,E,F$ denote the numbers of vertices, edges, and faces on the
sphere, then
\begin{equation}\label{eq:euler}
 3V=2E,\qquad V-E+F=2,\qquad V=2(F-2),\quad E=3(F-2).
\end{equation}
\end{proof}

\begin{lemma}[Area variations and pressures]\label{lem:area-surjectivity}
For a regular cluster with three finite cells of positive area in
$\mathbb H^2$, the area differential on compactly supported smooth
ambient vector fields is surjective onto $\mathbb R^3$.  At a fixed-area minimizer there is
a unique vector $p\in\mathbb R^3$ satisfying
\eqref{eq:first-variation}.
\end{lemma}

\begin{proof}
The label-adjacency graph has vertices $0,1,2,3$, with $i,j$
adjacent when
$\mathcal H^1(\partial^*E_i\cap\partial^*E_j)>0$.
Its vertices are labels rather than faces. It is connected.
Otherwise, the union of
the cells in one graph component would have zero distributional
gradient on the connected ambient surface, although both that union
and its complement have positive area.  Choose a spanning tree
on the vertices $0,1,2,3$ and choose a small regular interface
interval for each of its three edges, with disjoint supports.
A normal bump field on the oriented edge $i/j$ has area vector
$\gamma(e_i-e_j)$, where $\gamma\ne0$ and $e_0=0$.
The response matrix is
\[
 D\,\operatorname{diag}(\gamma_1,\gamma_2,\gamma_3),
\]
where $D$ is the tree incidence matrix with row $0$ removed.
Successively removing leaves different from $0$ shows $|\det D|=1$.
Thus the response matrix is invertible.  The corresponding three
flows give exact local area corrections by the inverse function theorem.

If a smooth variation has zero first variation of all three finite-cell areas, this correction
changes its parameter by $O(t^2)$ and makes the areas exactly
constant.  Minimality implies that its first perimeter derivative
vanishes.  Hence $\delta P$ vanishes on the kernel of $\delta v$
and factors uniquely as $p\cdot\delta v$.  Varying a single
regular interface interval gives $k_{ij}=p_i-p_j$ with the stated
outward-normal convention.
\end{proof}

\begin{lemma}[Gauss--Bonnet for a disk face]\label{lem:gauss-bonnet}
Let $D$ be a disk face of a regular equilibrium network.  If it
has $d$ sides, pressure $p_D$, and neighboring label $j(e)$ across
side $e$, then
\begin{equation}\label{eq:gauss-bonnet}
 \sum_{e\subset\partial D}(p_D-p_{j(e)})\ell(e)
 =2\pi+|D|-\frac{d\pi}{3}.
\end{equation}
\end{lemma}

\begin{proof}
Orient $\partial D$ with $D$ to the left.  If $n$ is the left
normal, then $n=-N$ and $\nabla^H_TT=-kN=kn$; the geodesic curvature
in Gauss--Bonnet is therefore $k=p_D-p_{j(e)}$.
The ambient Gaussian curvature is $-1$ and each interior angle
is $2\pi/3$.  Thus
\[
 -|D|+\sum_{e\subset\partial D}k_e\ell(e)
       +d\left(\pi-\frac{2\pi}{3}\right)=2\pi,
\]
which is \eqref{eq:gauss-bonnet}.
\end{proof}

\section{The standard hyperbolic family}
\label{sec:standard}
Fix $n\ge2$, write $m=n-1$, and let $1\le k\le n+1$.
We use the hyperboloid model
\[
 \mathbb H^n=\{y=(y^{\rm sp},y_0)\in\mathbb R^{n,1}:
       \langle y,y\rangle_L=-1,\ y_0>0\},\qquad
 \langle(x,t),(x',t')\rangle_L=x\cdot x'-tt'.
\]
We write $D$ for the flat ambient derivative and $\nabla^H$ for
the Levi--Civita covariant derivative of the hyperbolic metric.
A timelike vector has negative Lorentz squared norm; it is future
directed when its last coordinate is positive and past directed when
that coordinate is negative. The shape operator with unit normal
$\nu$ is $u\mapsto\nabla^H_u\nu$, and our mean-curvature
normalization is $H=\operatorname{tr}(\nabla^H\nu)$.
Thus a sphere with outward normal has mean curvature $m\coth R$.

Here is the standard family in arbitrary dimension. Choose unit
vectors $a_0,\ldots,a_k\in\mathbb R^{n+1}$ with
\[
 \sum_{i=0}^ka_i=0,\qquad a_i\cdot a_j=-\frac1k\quad(i\ne j),
 \qquad U_i=\{x\in\mathbb S^n:a_i\cdot x<a_j\cdot x
                                  \text{ for all }j\ne i\}.
\]
These are the cells of the regular-simplex spherical partition.
A common positive rescaling of the sites does not change these
inequalities; in the proofs we normalize their pairwise distances
to one.
The map $y\mapsto(y^{\rm sp}/y_0,1/y_0)$ identifies the hyperboloid
with the upper hemisphere
$\mathbb S^n_+=\{(x^{\rm sp},x_0)\in\mathbb S^n:x_0>0\}$,
equipped with $x_0^{-2}g_{\mathbb S^n}$.
Indeed, putting $x^{\rm sp}=y^{\rm sp}/y_0$ and $x_0=1/y_0$,
the identity $y^{\rm sp}\cdot dy^{\rm sp}=y_0\,dy_0$ gives
\[
 |dx|^2=\left|\frac{dy^{\rm sp}}{y_0}
            -\frac{y^{\rm sp}}{y_0^2}\,dy_0\right|^2
            +\frac{dy_0^2}{y_0^4}
        =\frac{|dy^{\rm sp}|^2-dy_0^2}{y_0^2}.
\]
A standard hyperbolic $k$-bubble is
\[
 \Omega_i=\Phi(U_i)\cap\mathbb S^n_+\quad(0\le i\le k),
 \qquad \overline{\Phi(U_i)}\subset\mathbb S^n_+
                      \quad(1\le i\le k),
\]
where $\Phi$ is a M\"obius transformation of $\mathbb S^n$.
The exterior label is $0$. For $R>0$, set
\[
 W_n(R)=|\mathbb S^{n-1}|\int_0^R\sinh^{n-1}t\,dt,
 \qquad A_n(R)=|\mathbb S^{n-1}|\sinh^{n-1}R.
\]
The perimeter of a ball of volume $a$ is
\[
 B_n(a)=A_n(W_n^{-1}(a)),\qquad B_n(0)=0.
\]
For $a>0$,
\[
 \begin{aligned}
 B_n'(a)&=\frac{A_n'(R)}{W_n'(R)}\\
 &=\frac{m|\mathbb S^{n-1}|\sinh^{m-1}R\cosh R}
        {|\mathbb S^{n-1}|\sinh^mR}=m\coth R.
 \end{aligned}
\]
Here $R=W_n^{-1}(a)>0$.
We write $T_k(v)$ for the perimeter of a standard cluster with finite-cell
volume vector $v$. The next results establish that this notation is well defined.

\begin{lemma}[Pressure realization]\label{lem:standard-pressure}
Every $\kappa\in(1,\infty)^k$ determines, up to a label-preserving
hyperbolic isometry, a unique standard hyperbolic $k$-bubble with
pressures $m\kappa_i$ in the finite cells and exterior pressure zero. Every bounded standard
$k$-bubble with positive finite-cell volumes has such a representation.
\end{lemma}
\begin{proof}
Let $\mathbf1=(1,\ldots,1)^T\in\mathbb R^k$ and define
\begin{equation}\label{eq:model-gram}
 \mathsf A=\frac{I_k+\mathbf1\mathbf1^T}{2},\qquad
 G=\mathsf A-\kappa\kappa^T,\qquad
 \alpha=\kappa^T\mathsf A^{-1}\kappa
 =2\sum_i\kappa_i^2-\frac{2}{k+1}\Bigl(\sum_i\kappa_i\Bigr)^2.
\end{equation}
Put $J_0=\mathbf1\mathbf1^T$, so $J_0^2=kJ_0$. Then
\[
 \mathsf A\,2\left(I_k-\frac{J_0}{k+1}\right)
 =(I_k+J_0)\left(I_k-\frac{J_0}{k+1}\right)
 =I_k+\left(1-\frac1{k+1}-\frac{k}{k+1}\right)J_0=I_k.
\]
Thus $\mathsf A^{-1}=2(I_k-J_0/(k+1))$, giving the formula for $\alpha$.
Cauchy--Schwarz gives
\[
 \alpha\ge\frac{2}{k+1}\sum_i\kappa_i^2
       >\frac{2k}{k+1}\ge1.
\]
If $u=\mathsf A^{-1/2}\kappa$, then
$\mathsf A^{-1/2}G\mathsf A^{-1/2}=I_k-uu^T$, whose eigenvalues are $1-\alpha<0$
and $1$ with multiplicity $k-1$. Thus $G$ is nondegenerate with signature
$(k-1,1)$. Since $k-1\le n$, it has a realization by independent vectors
$c_1,\ldots,c_k\in\mathbb R^{n,1}$ satisfying
\begin{equation}\label{eq:model-centers}
 \langle c_i,c_j\rangle_L=\mathsf A_{ij}-\kappa_i\kappa_j.
\end{equation}
The diagonal entries are $1-\kappa_i^2<0$, and the off-diagonal entries
are $1/2-\kappa_i\kappa_j<0$. Consequently all the $c_i$ are timelike
with the same time orientation; choose the past orientation. An isometry
between two such spans extends by an orthogonal isometry of their
positive-definite complements. The extension preserves the future cone,
so it induces a hyperbolic isometry. This proves uniqueness of the realization.

Set $c_0=0$, $\kappa_0=0$, and define
\begin{equation}\label{eq:model-scores}
 F_0=0,\qquad F_i(y)=\langle c_i,y\rangle_L-\kappa_i,
 \qquad \Omega_i=\{F_i<F_j\text{ for all }j\ne i\}.
\end{equation}
For distinct labels, write $F_{ij}=F_i-F_j$ and
\[
 \Sigma_{ij}=\{y\in\mathbb H^n:F_i(y)=F_j(y)<F_\ell(y)
                         \text{ for all }\ell\notin\{i,j\}\}.
\]
This is the relatively open $(n-1)$-dimensional interface; integrals
over it use $dA=d\mathcal H^{n-1}$. Intersections of three or more
cells are not included in $\Sigma_{ij}$.
Distinct costs do not agree on an open set, so their tie sets have zero
$n$-dimensional volume and these cells form a partition up to null sets.
On an interface $F_i=F_j$, its outward normal from $\Omega_i$ is
\[
 \nu_{ij}=c_i-c_j+(\kappa_i-\kappa_j)y.
\]
Write $d=c_i-c_j$ and $\delta=\kappa_i-\kappa_j$. On the interface
$\langle d,y\rangle_L=\delta$, so
\[
 \langle d+\delta y,y\rangle_L=\delta-\delta=0,\qquad
 |d+\delta y|_L^2=\langle d,d\rangle_L+2\delta^2-\delta^2
                 =\langle d,d\rangle_L+\delta^2=1.
\]
The sign is outward because $\Omega_i$ is the side $F_i-F_j<0$.
For a vector $u$ tangent to the interface, ambient differentiation gives
$D_u\nu_{ij}=\delta u$. This vector is already tangent to
$\mathbb H^n$, so $\nabla^H_u\nu_{ij}=\delta u$, giving the shape
operator $\delta I$.

To identify the partition, use the conformal compactification
\[
 p=(y^{\rm sp}/y_0,1/y_0)\in\mathbb S^n_+,
 \qquad \beta_i=(c_i^{\rm sp},-\kappa_i,c_{i0})\in\mathbb R^{n+1,1}.
\]
With $\beta_0=0$, the full-sphere costs are
\[
 \ell_i(p)=\langle \beta_i,(p,1)\rangle_{n+1,1}
       =c_i^{\rm sp}\cdot p^{\rm sp}-\kappa_i p_0-c_{i0},
 \qquad \ell_i(p)=F_i(y)/y_0\quad(p_0>0).
\]
Their Gram matrix is $\mathsf A$, hence the span of $\beta_1,\ldots,\beta_k$ is
positive definite. Choose a future-preserving Lorentz transformation
that takes this span into a spatial hyperplane. Under its projective
action on future null rays, the sphere is mapped conformally to itself:
if $\Lambda(p,1)=\lambda(p)(p',1)$, then $\lambda(p)>0$, and Lorentz
invariance gives
\[
 \ell_i(p)=\lambda(p)\langle\Lambda \beta_i,(p',1)\rangle_{n+1,1}.
\]
Thus the minimizing-cost inequalities are preserved. The transformed
sites $0,\Lambda \beta_1,\ldots,\Lambda \beta_k$ are spatial and have every
pairwise squared distance equal to one. After a common translation of
these sites, which does not change any comparison of costs, their
partition is the regular-simplex partition. This identifies the full
spherical partition as a standard partition. In particular every cell
and every pair interface is nonempty, and every interface has positive
$(n-1)$-dimensional measure.

Past timelikeness gives $-c_{i0}>|c_i^{\rm sp}|$. On $p_0\le0$,
\[
 \ell_i(p)\ge-|c_i^{\rm sp}|-c_{i0}>0=\ell_0(p).
\]
By compactness and continuity, this strict inequality holds on a whole
neighborhood of the closed southern hemisphere. Every finite cell and
its closure is consequently compactly contained in $\mathbb S^n_+$.
The resulting hyperbolic cluster is bounded.

Conversely, let a standard spherical partition have all its finite
cells compactly contained in $\mathbb S^n_+$. Start with spatial
regular-simplex sites of common pairwise distance one and apply a
Lorentz lift of its defining M\"obius transformation. Subtract the
transformed exterior site from all sites. Their differences $\beta_i$ then
satisfy $\langle \beta_i,\beta_j\rangle_{n+1,1}=\mathsf A_{ij}$. Write
$\beta_i=(c_i^{\rm sp},-\kappa_i,c_{i0})$; dividing the sphere costs by the
positive factor $1/y_0$ gives exactly \eqref{eq:model-scores} and
\eqref{eq:model-centers}.

Let $D_i=\{F_i<0\}$. At every point with a unique minimizing cost,
\begin{equation}\label{eq:model-support-balls}
 \bigcup_{i=1}^k\Omega_i=\bigcup_{i=1}^kD_i,\qquad
 \Omega_i\subset D_i\subset\bigcup_{j=1}^k\Omega_j.
\end{equation}
The omitted tie sets are null. Since the finite cells are positive and
bounded, each open set $D_i$ is nonempty and bounded. We check that this
forces $c_i$ to be past timelike and $\kappa_i>1$. Write $c_i=(a,t)$.
If $t>-|a|$, choose a unit $\theta$ with $a\cdot\theta=-|a|$; along
$y(R)=(\sinh R\,\theta,\cosh R)$, the cost tends to $-\infty$.
Thus $D_i$ would be unbounded. If $t=-|a|$ and $c_i\ne0$, past
nullness gives $\langle c_i,y\rangle_L>0$ on $\mathbb H^n$. For
$\kappa_i\le0$ the set is empty, whereas for $\kappa_i>0$ the same
ray has cost $|a|e^{-R}-\kappa_i<0$ for large $R$, again giving an
unbounded set. If $c_i=0$, the set is either empty or all of
$\mathbb H^n$. These alternatives exclude every case except
$t<-|a|$. A hyperbolic isometry now takes $c_i$ to $(0,-a_0)$ with
$a_0>0$. Since $a_0^2=\kappa_i^2-1$, nonemptiness of
\[
 D_i=\{a_0\cosh R<\kappa_i\}
\]
implies $\kappa_i>a_0>0$, hence $\kappa_i>1$. Its radius satisfies
\[
 \cosh R_i=\frac{\kappa_i}{\sqrt{\kappa_i^2-1}},\qquad
 \sinh R_i=\frac1{\sqrt{\kappa_i^2-1}},\qquad
 \coth R_i=\kappa_i.
\]
Thus $D_i$ is the ball of radius
$R_i=\operatorname{arccoth}\kappa_i$. The preceding
Gram realization proves the asserted uniqueness.
\end{proof}

\begin{theorem}[Global volume coordinates]\label{thm:standard-volume}
The map $\kappa\mapsto v(\kappa)$ is a smooth diffeomorphism from
$(1,\infty)^k$ onto $(0,\infty)^k$. There are positive weights
$w_{ij}=w_{ji}$. Let $L_w$ be the grounded weighted Laplacian of the
complete label graph on $0,\ldots,k$, obtained by deleting the
exterior row and column. With $b_0=0$, its formula is
\[
 (L_wb)_i=\sum_{\substack{0\le j\le k\\j\ne i}}w_{ij}(b_i-b_j),
 \qquad
 D_\kappa v=-L_w\prec0,\quad
 D_vT_k=m\kappa,\quad D_v^2T_k=-mL_w^{-1}\prec0.
\]
Writing $s=\sum_i v_i$, one has
\begin{equation}\label{eq:model-derivative-bounds}
 B_n'(s)\le\partial_iT_k(v)\le B_n'(v_i),\qquad
 \partial_iT_k(v)>m,\qquad ds=-\sum_iw_{i0}\,d\kappa_i.
\end{equation}
In particular, decreasing pressures coordinatewise, with at least one
strict decrease, strictly increases the total finite volume.
\end{theorem}
\begin{proof}
Let $C$ have columns $c_i$ and let $J_L$ be the matrix of the Lorentz
form, so $C^TJ_LC=G$. Define
\[
 t=CG^{-1}\kappa,\qquad V(y)=-\langle t,y\rangle_L-1.
\]
The rank-one inverse formula gives
\[
 G^{-1}=\mathsf A^{-1}+\frac{\mathsf A^{-1}\kappa\kappa^T\mathsf A^{-1}}{1-\alpha},
 \qquad
 C^TJ_Lt=(C^TJ_LC)G^{-1}\kappa=GG^{-1}\kappa=\kappa.
\]
Moreover
\[
 \begin{split}
 \langle t,t\rangle_L
 &=\kappa^TG^{-1}(C^TJ_LC)G^{-1}\kappa
   =\kappa^TG^{-1}\kappa\\
 &=\alpha+\frac{\alpha^2}{1-\alpha}
   =\frac{\alpha(1-\alpha)+\alpha^2}{1-\alpha}
   =-\frac{\alpha}{\alpha-1}.
 \end{split}
\]
Thus $t$ is timelike. The signs $\langle c_i,t\rangle_L=\kappa_i>0$
and the past orientation of $c_i$ show that $t$ is future directed.
For future unit timelike $y$, the reverse Lorentz Cauchy--Schwarz
inequality yields
\begin{equation}\label{eq:model-positive-weight}
 V(y)\ge\sqrt{\frac{\alpha}{\alpha-1}}-1>0.
\end{equation}

We first justify smooth parameter dependence. Near any $\kappa^*$,
put $\kappa_\tau=\kappa^*+\tau(\kappa-\kappa^*)$ and solve
\[
 \frac{dC_\tau}{d\tau}
   =-C_\tau G(\kappa_\tau)^{-1}\kappa_\tau(\kappa-\kappa^*)^T,
 \qquad C_0=C^*.
\]
The coefficients are smooth on a compact parameter neighborhood, so
the solution is smooth in $(\tau,\kappa)$. For a general pressure path
with derivative $b$, the error $E=C^TJ_LC-G$ satisfies
\[
 E'=-b\kappa^TG^{-1}E-EG^{-1}\kappa b^T,
 \qquad E(0)=0.
\]
Uniqueness for this linear matrix equation gives $E=0$. The endpoint
therefore supplies a smooth local choice of centers $C(\kappa)$.

The lifted vectors $\beta_i(\kappa)$ have the fixed positive Gram matrix
$\mathsf A$. Complete an orthonormal basis of their span to a Lorentz
orthonormal frame, locally smoothly in $\kappa$: project a fixed
future timelike vector to the orthogonal complement and normalize it,
then apply ordinary Gram--Schmidt in the remaining positive-definite
complement. The resulting smooth Lorentz transformations carry the
lifted sites at $\kappa^*$ to those at $\kappa$. Their projective
actions give smooth sphere diffeomorphisms taking each fixed labeled
cell and interface to its counterpart. By the strict southern-hemisphere
inequality in Lemma~\ref{lem:standard-pressure}, the closures of all
finite cells in a sufficiently small parameter neighborhood lie in a
single compact subset of $\mathbb S^n_+$. On that compact set the
hyperbolic density and all differentiated Jacobians are smooth and
bounded. Pullback to the fixed cells and interfaces thus permits
differentiation of volume and perimeter to every finite order.

Fix a pressure direction $b\in\mathbb R^k$ and set $b_0=0$.
For this fixed direction, solve the pressure-path ODE afresh with
$\kappa'=b$. Uniqueness up to isometry in
Lemma~\ref{lem:standard-pressure} ensures that it computes the same
volume map $v(\kappa)$. This path realizes the infinitesimal motion
$\dot c_i=-b_i t$. On a moving interface,
\[
 \dot F_{ij}=(b_i-b_j)V,\qquad
 |\nabla F_{ij}|=1,
 \qquad f_{ij}=-(b_i-b_j)V,
\]
where $f_{ij}$ is its outward normal speed from cell $i$.
The volume transport formula gives
\[
 (D_\kappa v\,b)_i=-\sum_{j\ne i}w_{ij}(b_i-b_j),
 \qquad w_{ij}=\int_{\Sigma_{ij}}V\,dA.
\]
Each weight is finite and strictly positive by boundedness,
\eqref{eq:model-positive-weight}, and positivity of each interface area.
Consequently
\[
 b^TL_wb=\sum_{0\le i<j\le k}w_{ij}(b_i-b_j)^2>0
 \quad\text{for }b\ne0,\ b_0=0.
\]
The last strict inequality follows already from the terms $w_{i0}b_i^2$.
The inverse function theorem now gives local invertibility.

The supporting ball $D_i$ from \eqref{eq:model-support-balls} has radius
$R_i=\operatorname{arccoth}\kappa_i$, and the inclusions there give
\begin{equation}\label{eq:model-ball-volume}
 v_i\le W_n(R_i)\le\sum_jv_j=s.
\end{equation}
If $\delta\le v_i\le M$ for all $i$, then
\[
 \coth(W_n^{-1}(kM))\le\kappa_i
       \le\coth(W_n^{-1}(\delta)).
\]
Thus the inverse image of every compact subset of $(0,\infty)^k$ is
compact in $(1,\infty)^k$. The map is proper. For distinct
$\kappa,\lambda$, let $h=\kappa-\lambda$. The whole segment is in
the pressure domain, and
\[
 h\cdot\bigl(v(\kappa)-v(\lambda)\bigr)
 =-\int_0^1h^TL_w(\lambda+\tau h)h\,d\tau<0.
\]
Hence the map is injective. Its image is nonempty and open by local
invertibility. Properness makes it relatively closed in the positive
volume orthant: a convergent image sequence in that orthant has a
convergent preimage subsequence. The orthant is connected, so the image
is the whole orthant. The local smooth inverses combine to a global
smooth inverse.

The regular interfaces have mean curvature $m(\kappa_i-\kappa_j)$,
and their conormals cancel along every triple stratum because the
angles are $120$ degrees. Higher strata have zero
$(n-2)$-dimensional measure. Summing the first area variation therefore
gives
\[
 dT_k=m\sum_i\kappa_i\,dv_i.
\]
One may apply this formula to the preceding parameter variations by
extending their ambient velocities with a cutoff equal to one on the
common compact closure of the finite cluster. Since
$D_\kappa v=-L_w$, differentiation in volume
coordinates gives $D_vT_k=m\kappa$ and $D_v^2T_k=-mL_w^{-1}$.
Apply the decreasing function $B_n'$ to
\eqref{eq:model-ball-volume} to obtain
\eqref{eq:model-derivative-bounds}. Summing the volume variations gives
\[
 \begin{split}
 ds&=-\sum_{1\le i<j\le k}w_{ij}
       \bigl[(d\kappa_i-d\kappa_j)+(d\kappa_j-d\kappa_i)\bigr]
       -\sum_iw_{i0}\,d\kappa_i\\
   &=-\sum_iw_{i0}\,d\kappa_i.
 \end{split}
\]
Its integration along a coordinatewise
decreasing pressure segment proves strict total-volume monotonicity.
\end{proof}

\begin{proposition}[Coordinate-face limits]\label{prop:standard-faces}
The function $T_k$ has a continuous extension to $[0,\infty)^k$.
At a vector with exactly $r$ positive coordinates its value is the
standard $r$-bubble perimeter of those coordinates. Set $T_0=0$.
Then, for $u\in[0,\infty)^{k-1}$ and $a\ge0$,
\begin{equation}\label{eq:model-face-bound}
 0\le T_k(u,a)-T_{k-1}(u)\le B_n(a).
\end{equation}
\end{proposition}
\begin{proof}
For positive $v,w$, integrate the coordinate derivative bounds along
the $k$ segments changing one coordinate at a time. This gives
\begin{equation}\label{eq:model-continuity-modulus}
 |T_k(v)-T_k(w)|\le\sum_{i=1}^k|B_n(v_i)-B_n(w_i)|.
\end{equation}
Since $B_n$ is continuous on $[0,\infty)$, every positive sequence
converging to a finite nonnegative vector has Cauchy perimeter values,
and any two such sequences have the same limit. This defines the
unique continuous extension and preserves the displayed modulus.

We identify its values. Fix a standard $r$-bubble with pressure vector
$\bar\kappa$ and centers $\bar c_a$, $1\le a\le r$. For $L$ large,
construct the standard $k$-bubble with the first $r$ pressures
$\bar\kappa_a$ and the remaining $k-r$ pressures equal to $L$.
If $r\ge1$, the Gram matrix of its first $r$ centers is the fixed
nondegenerate matrix of the $\bar c_a$. A hyperbolic isometry aligns
these centers with the fixed centers. Removing the additional scores
then recovers the original $r$-bubble exactly. Each extra label has
a supporting ball of radius $R_L=\operatorname{arccoth}L$, with
$R_L\to0$. Each extra-cell volume is at most $W_n(R_L)$, and each
surviving cell loses volume at most $(k-r)W_n(R_L)$. The volume vectors
therefore converge to the desired coordinate-face vector.

There are only finitely many label pairs. An interface between an
extra label and a surviving label $a$ (including the exterior,
$\bar\kappa_0=0$) has absolute principal curvature $L-\bar\kappa_a>1$.
Its entire supporting sphere has area
\[
 A_n\bigl(\operatorname{arccoth}(L-\bar\kappa_a)\bigr)\longrightarrow0.
\]
An interface between two extra labels is totally geodesic and lies
in the closure of an extra supporting ball. The intersection of a
totally geodesic hyperplane with an ambient radius-$R_L$ ball is either
empty or an intrinsic ball of radius at most $R_L$: if the center has
distance $d$ from the hyperplane, its orthogonal projection satisfies
$\cosh\operatorname{dist}(c,x)=\cosh d\,
\cosh\operatorname{dist}(\pi(c),x)$. Hence this interface area is at
most
\[
 |\mathbb S^{n-2}|\int_0^{R_L}\sinh^{n-2}t\,dt\longrightarrow0.
\]
Every surviving--surviving interface is a subset of the corresponding
fixed interface. The removed part lies in the union of the closed extra
supporting balls. Each fixed interface has compact closure in a smooth
supporting hypersurface. A finite covering by graphical coordinate
patches gives constants $C,r_*>0$ such that its area in every ambient
ball of radius $r<r_*$ is at most $Cr^{n-1}$, uniformly in the center
of the ball. Thus the total removed area is $O(R_L^{n-1})$.
All new and removed interface areas tend to zero, proving perimeter
convergence to the fixed $r$-bubble perimeter. For $r=0$, all labels
are extra and the same estimates give convergence to zero. By the
unique continuous extension, these sequences identify every face value.

For positive coordinates, integrate
$0<\partial_kT_k\le B_n'(v_k)$ from $v_k=\varepsilon$ to $a$,
then let $\varepsilon\downarrow0$ and use the identified face value.
Approximation by positive values of the other coordinates gives
\eqref{eq:model-face-bound} in its full stated range.
\end{proof}

\begin{corollary}[Standard excess]\label{cor:standard-excess}
Set $j_n(a)=B_n(a)-ma$ and
$S_{n,k}(v)=T_k(v)-m\sum_i v_i$. The function $S_{n,k}$ is symmetric,
smooth and strictly concave on $(0,\infty)^k$, extends continuously to
the coordinate faces, and satisfies
\[
 \begin{gathered}
 0<\partial_iS_{n,k}(v)=m(\kappa_i-1)\le j_n'(v_i),\\
 0\le S_{n,k}(u,a)-S_{n,k-1}(u)\le j_n(a),\\
 0<S_{n,k}(v)\le\sum_i j_n(v_i).
 \end{gathered}
\]
For $n=2$, write $S_k=S_{2,k}$ and $j=j_2$. Then
\[
 j(a)=\sqrt{a^2+4\pi a}-a
     =\frac{4\pi a}{\sqrt{a^2+4\pi a}+a}<2\pi,
 \qquad S_k(v)<2\pi k.
\]
\end{corollary}
\begin{proof}
Relabeling the finite cells proves symmetry. Subtracting a linear
function preserves the negative-definite Hessian. The derivative
claims follow from Theorem~\ref{thm:standard-volume}, and integration
in the last coordinate followed by the face limit proves the increment
bound. Starting at the origin and adding the coordinates successively
gives $0\le S_{n,k}(v)\le\sum_i j_n(v_i)$. On the positive orthant,
all coordinate derivatives are strictly positive; in particular
$S_{n,k}(v)>S_{n,k}(v/2)\ge0$. The displayed formula for $j$ and the
strict inequality $\sqrt{a^2+4\pi a}+a>2a$ give the planar bounds.
\end{proof}
\begin{lemma}[Hyperbolic isoperimetric inequality]\label{lem:hyperbolic-isoperimetric}
Let $n\ge2$ and let $E\subset\mathbb H^n$ have finite volume and
finite perimeter. Then
\[
 \operatorname{Per}_H(E)\ge B_n(|E|_H).
\]
\end{lemma}
\begin{proof}
For a compact set with Lipschitz boundary, this is the classical
space-form inequality in
\cite[Theorem 4.5.1(c)]{BFR}. We pass to finite-perimeter sets.
The case $|E|_H=0$ is immediate; assume $|E|_H>0$.
Fix a center and put $E_R=E\cap B_R$. The coarea identity gives
\[
 \int_0^\infty
 \mathcal H_H^{n-1}(E^{(1)}\cap\partial B_R)\,dR=|E|_H<\infty.
\]
Here $E^{(1)}$ is the set of points of density one. Choose
$R_h\uparrow\infty$ outside the exceptional radii for BV slicing, with
$\mathcal H_H^{n-1}(E^{(1)}\cap\partial B_{R_h})\to0$.
The intersection perimeter estimate and lower semicontinuity yield
\[
 \begin{split}
 \operatorname{Per}_H(E)
 &\le\liminf_h\operatorname{Per}_H(E_{R_h})\\
 &\le\limsup_h\left[
     \operatorname{Per}_H(E;B_{R_h})
     +\mathcal H_H^{n-1}(E^{(1)}\cap\partial B_{R_h})\right]
 \le\operatorname{Per}_H(E).
 \end{split}
\]
Also $|E_{R_h}|_H\to|E|_H$ by monotone convergence. It remains
to approximate a bounded finite-perimeter set $F$.

Choose a Poincar\'e coordinate ball containing $\overline F$ in
its interior. On a slightly larger compact coordinate set $K$,
$g_H=\rho^2g_E$ with $\rho>0$ smooth and bounded, and
\[
 |F|_H=\int_F\rho^n\,dx,\qquad
 \operatorname{Per}_H(F)
       =\int_K\rho^{n-1}\,d|D\mathbf1_F|.
\]
Smooth approximation of finite-perimeter sets
\cite[Theorem 13.8 and Remark 13.9]{Maggi}, using convolution
and regular level sets in this coordinate ball, gives smooth
bounded sets $F_j\Subset K$ such that
\[
 \mathbf1_{F_j}\longrightarrow\mathbf1_F\quad\text{in }L^1(dx),
 \qquad |D\mathbf1_{F_j}|\stackrel{*}{\rightharpoonup}|D\mathbf1_F|.
\]
The approximations may have support in this fixed $K$ because the
convolutions of $\mathbf1_F$ have support in shrinking neighborhoods
of $\overline F$. Applying these convergences to the bounded
continuous weights $\rho^n$ and $\rho^{n-1}$ gives
\[
 |F_j|_H\to|F|_H,\qquad
 \operatorname{Per}_H(F_j)\to\operatorname{Per}_H(F).
\]
The closures $\overline{F_j}$ satisfy the compact Lipschitz-boundary
inequality, and their boundaries have zero $n$-volume. Therefore
\[
 \operatorname{Per}_H(F_j)\ge B_n(|F_j|_H).
\]
Letting $j\to\infty$, then applying this result to $E_{R_h}$ and
letting $h\to\infty$, proves the assertion because $B_n$ is continuous.
\end{proof}

\section{Circular networks: rigidity and local deletion}
\label{sec:circular-networks}

We use the network terminology of Section~\ref{sec:preliminaries}.
Faces and their boundary cycles are viewed on the one-point
compactification when the exterior face is included. A generalized
Euclidean circle means a circle or a straight line. The local
deletion statements apply to equilibrium networks; the rolling
argument also uses fixed-area minimality.

\begin{lemma}[$K_4$ rigidity]\label{lem:k4}
Let a bounded regular equilibrium three-bubble in $\mathbb H^2$ have
four complementary disk faces, four vertices and six edges, with
embedded graph $K_4$. Then it is standard. The same conclusion holds
for a bounded Euclidean circular-arc equilibrium network with this
topology.
\end{lemma}
\begin{proof}
Use the hyperboloid model with spatial coordinates first and Lorentz
product $\langle(x,t),(x',t')\rangle_L=x\cdot x'-tt'$.
For an oriented $i/j$ edge let $N_{ij}$ point from label $i$ into
label $j$ and put $\kappa_{ij}=p_i-p_j$. Its ambient quasi-center
\[
 c_{ij}=N_{ij}-\kappa_{ij}y
\]
is constant. Indeed, differentiation of
$\langle N_{ij},y\rangle_L=0$ gives
$\langle D_TN_{ij},y\rangle_L=-\langle N_{ij},T\rangle_L=0$.
Thus $D_TN_{ij}$ is tangent to the hyperboloid and equals
$\nabla^H_TN_{ij}=\kappa_{ij}T$. Since $D_Ty=T$ and
$\kappa_{ij}$ is constant along the edge,
$D_Tc_{ij}=\kappa_{ij}T-\kappa_{ij}T=0$.
At a triple vertex, cyclic normals and pressure differences sum to
zero, and hence
\[
 c_{ij}+c_{j\ell}+c_{\ell i}=0.
\]
Each of the four labels occupies one face, and all four dual
triangles occur. Taking label $0$ as reference gives vectors $c_i$
and numbers $\kappa_i=p_i-p_0$ with
\[
 c_0=0,\quad \kappa_0=0,\quad
 c_{ij}=c_i-c_j,\quad
 \langle c_i-c_j,c_i-c_j\rangle_L+(\kappa_i-\kappa_j)^2=1.
\]
The last identity follows by expanding
\[
 \langle N_{ij}-\kappa_{ij}y,N_{ij}-\kappa_{ij}y\rangle_L
 =1-2\kappa_{ij}\langle N_{ij},y\rangle_L
       +\kappa_{ij}^2\langle y,y\rangle_L
 =1-\kappa_{ij}^2.
\]

Write $c_i=(c_i^{\rm sp},c_i^{\rm tm})$. In the conformal lift of
Section~\ref{sec:standard}, the vectors
\[
 b_i=(c_i^{\rm sp},-\kappa_i,c_i^{\rm tm})\in\mathbb R^{3,1},
 \qquad b_0=0,
\]
satisfy $\langle b_i-b_j,b_i-b_j\rangle_{3,1}=1$ for every pair.
To see both the equation of the lifted support and its norm, put
\[
 p=\left(\frac{y^{\rm sp}}{y^{\rm tm}},\frac1{y^{\rm tm}}\right)
 \in\mathbb S^2,
 \qquad z=(p,1).
\]
The hyperboloid identity gives
$|p|^2=(|y^{\rm sp}|^2+1)/(y^{\rm tm})^2=1$, and
\[
 \langle b_i-b_j,z\rangle_{3,1}
 =\frac{\langle c_i-c_j,y\rangle_L-(\kappa_i-\kappa_j)}
        {y^{\rm tm}},\qquad
 \langle b_i-b_j,b_i-b_j\rangle_{3,1}
 =\langle c_i-c_j,c_i-c_j\rangle_L+(\kappa_i-\kappa_j)^2.
\]
Since $b_0=0$, each $b_i$ has norm squared one, and polarization gives
$\langle b_i,b_j\rangle_{3,1}=(1+1-1)/2=1/2$ for $i\ne j$.
The three nonzero vectors consequently have Gram matrix
\[
 \bigl(\langle b_i,b_j\rangle_{3,1}\bigr)_{i,j=1}^3
 =\frac12(I_3+\mathbf1\mathbf1^T)\succ0.
\]
For $a\ne0$, its quadratic form is
$\tfrac12|a|^2+\tfrac12(\sum_i a_i)^2>0$.
Their span therefore has a future-pointing unit timelike normal.
An orthonormal basis of the span, completed by that normal, gives
a time-orientation-preserving Lorentz transformation carrying the
span to a spatial hyperplane. Its action on the conformal sphere
is a global M\"obius diffeomorphism, and all six supporting circles
become great circles.

This transformation applies to the entire embedded graph on the
sphere. The original finite faces have compact closure in the
hyperbolic domain, so extending label $0$ across its complementary
hemisphere adds no interface. The transformed graph therefore still
has exactly the four labeled faces and the six original edges.

Let the resulting four spatial sites be $w_i$; their pairwise
distances are one. With $\bar w=\tfrac14\sum_iw_i$, the variance
identity gives
\[
 \sum_i|w_i-\bar w|^2
 =\frac14\sum_{i<j}|w_i-w_j|^2=\frac32,
 \qquad |w_i-\bar w|^2=\frac34-\frac14\frac32=\frac38.
\]
Thus $(w_i-\bar w)\cdot(w_j-\bar w)=(3/8+3/8-1)/2=-1/8$
for $i\ne j$. Setting $v_i=\sqrt{8/3}(w_i-\bar w)$ gives
regular tetrahedral unit vectors $v_0,v_1,v_2,v_3$ with
\[
 \sum_i v_i=0,\qquad v_i\cdot v_j=-\frac13\quad(i\ne j).
\]
The differences $v_j-v_\ell$, with $j,\ell\ne i$, span the
two-dimensional plane $v_i^\perp$. The vertex at which label $i$ is absent is orthogonal to that plane, and is thus
$y_i=s_iv_i$, $s_i\in\{1,-1\}$.
For $j\ne i$ put
\[
 u_{ij}=\frac{v_j+v_i/3}{\sqrt{8/9}}.
\]
In the common tangent planes at $v_i$ and $-v_i$ these are unit
vectors: $|v_j+v_i/3|^2=1+1/9-2/9=8/9$. For distinct $j,\ell$
their product is
\[
 u_{ij}\cdot u_{i\ell}
 =\frac{-1/3-1/9-1/9+1/9}{8/9}=-\frac12.
\]
Since $(s_iv_i)\cdot(s_jv_j)=-s_is_j/3$, the short-arc tangent is
\[
 \frac{s_jv_j-[(s_iv_i)\cdot(s_jv_j)]s_iv_i}
      {\sqrt{1-[(s_iv_i)\cdot(s_jv_j)]^2}}
 =\frac{s_j(v_j+v_i/3)}{\sqrt{8/9}}=s_ju_{ij}.
\]
Let $\epsilon_{ij}=1$ for the short arc and $-1$ for its long
complement. Embeddedness makes each edge one of these two simple
arcs. The actual outgoing tangent is $\epsilon_{ij}s_ju_{ij}$.
The product of the tangents for edges $ij$ and $i\ell$ is
$-\tfrac12(\epsilon_{ij}s_j)(\epsilon_{i\ell}s_\ell)$.
For this to be $-1/2$, the three signs $\epsilon_{ij}s_j$ must
coincide, say with $t_i$. At the other endpoint
$t_j=\epsilon_{ij}s_i$, so $t_is_i=t_js_j$. Connectivity gives
\begin{equation}\label{eq:k4-signs}
 \epsilon_{ij}s_j=t_i,\qquad
 t_is_i=t_js_j=\sigma\in\{1,-1\},\qquad
 \epsilon_{ij}=\sigma s_is_j.
\end{equation}

At most one long edge can be incident to any vertex: two would
both contain its antipode in their interiors. That antipode is
different from all the other selected vertices.
Up to simultaneous sign reversal and permutation, the sign splits
are $4+0$, $3+1$ and $2+2$. With $a$ positive and $b$ negative
signs, $\sigma=1$ makes exactly the cross-sign edges long, of
degrees $b$ and $a$ in the two groups. For $\sigma=-1$ exactly
the same-sign edges are long, of degrees $a-1$ and $b-1$.
The maximal degrees for the three splits are consequently $0,3,2$
when $\sigma=1$, and $3,2,1$ when $\sigma=-1$.
Only the all-short case and the $2+2$ case with $\sigma=-1$ remain.
In the latter case choose $s_0=s_1=1$, $s_2=s_3=-1$. The short
arcs from $v_0$ to $v_1$ and from $-v_2$ to $-v_3$ have common
midpoint
\[
 q=\frac{v_0+v_1}{|v_0+v_1|},
\]
because $v_0+v_1=-(v_2+v_3)$. Their required long complementary
arcs both contain $-q$, contradicting embeddedness.

The remaining graph consists of all short arcs joining the $v_i$,
or its antipodal image. The first graph is precisely the boundary
of the partition
\[
 R_i=\{x\in\mathbb S^2:v_i\cdot x<v_j\cdot x
                         \text{ for every }j\ne i\}.
\]
For example, on the short edge between the vertices missing $i$
and $j$, a point is $(av_i+bv_j)/|av_i+bv_j|$ with $a,b>0$;
the other two scores equal $-(a+b)/(3|av_i+bv_j|)$, while the
scores of $i,j$ exceed that value by
$4a/(3|av_i+bv_j|)$ and $4b/(3|av_i+bv_j|)$, respectively.
Face labels therefore match this partition. The antipodal
case uses the sites $-v_i$. Undoing the M\"obius transformation
identifies the original labeled cluster as standard.

For the Euclidean assertion, the hypotheses give an equilibrium
triple bubble with six edges and four triple junctions, with precisely
the combinatorial type of the standard triple bubble. Lawlor's
rigidity statement \cite[Proposition 2.6]{Lawlor} therefore applies.
\end{proof}

\begin{lemma}[Digon deletion]\label{lem:digon-deletion}
Let $D$ be a bounded two-sided disk face of a regular equilibrium
network, with label $i$. Its two neighboring labels $j,k$ are
distinct. An arc $\gamma$ of signed curvature $p_j-p_k$ joins its
two vertices, has its relative interior in $D$, divides $D$ into
nonempty open pieces $D_j,D_k$, and matches the two external
$j/k$ edges smoothly. Replacing $D$ by these pieces produces an
embedded equilibrium network with the remaining pressures unchanged.
\end{lemma}
\begin{proof}
The three labels at either endpoint are distinct. Let the endpoints
be $A,B$. The quasi-centers $c_{ij},c_{ik}$ of the sides are
constant, and the two external edges both have quasi-center
\[
 c_{jk}=c_{ik}-c_{ij},\qquad
 p_j-p_k=(p_i-p_k)-(p_i-p_j).
\]
They consequently belong to the same complete generalized circle
$C$ through $A,B$, with the same oriented curvature.

On the compactified conformal plane, send $A$ to infinity and $B$
to zero by a M\"obius map. All three supporting circles become
lines through zero. The two sides of $D$ become rays, and the
selected domain becomes their sector of angle $2\pi/3$; this is
the interior angle at $B$. The external $j/k$ ray points away
from this sector, while its opposite ray is the interior bisector.
The inverse image of that opposite ray, together with its endpoints,
is the required arc $\gamma$. It lies in the closure of the old
face and divides it into two nonempty open domains.
The common supporting curve and the endpoint tangent data give
smooth matching at both ends. The inserted arc inherits its
oriented hyperbolic curvature from the unchanged supporting curve.
Every modification lies within $\overline D$.
\end{proof}

\begin{lemma}[Rolling excludes a digon]\label{lem:no-digon}
A fixed-area minimizing cluster with compact connected regular
boundary and at least four complementary faces has no bounded
digon.
\end{lemma}
\begin{proof}
Delete a supposed digon as in Lemma~\ref{lem:digon-deletion} and
write $N'$ for the resulting fixed network. The operation smooths
the two former vertices and joins their external $j/k$ edges.
It preserves connectedness and reduces the face count by one.
Thus $N'$ is connected and has at least three faces.
The maximal smooth edge containing $\gamma$ has a vertex endpoint:
a vertex-free closed component would, by connectedness, be the
whole network and have only two complementary faces.

Let $C$ be the complete supporting curve of $\gamma$. It has a
side-preserving one-parameter group of hyperbolic isometries $G_t$
acting on $C$ by arclength translation. The relation
$\langle c,c\rangle_L=1-\kappa^2$ from the quasi-center formula
classifies a complete support as a circle, a horocycle, or a
geodesic/equidistant curve according as this Lorentz norm is
negative, zero, or positive.
For a circle of radius $R$,
polar coordinates $d\rho^2+\sinh^2\rho\,d\theta^2$ give
$G_t(\rho,\theta)=(\rho,\theta+t/\sinh R)$, with unit speed
on $C$. A geodesic or equidistant curve can be placed at
$y_1=\sinh d$ in hyperboloid coordinates $(y_1,y_2,y_0)$.
Its parametrization and the preserving boost are
\[
 \begin{aligned}
 y(u)&=(\sinh d,\cosh d\sinh u,\cosh d\cosh u),\\
 B_\tau\begin{pmatrix}y_1\\y_2\\y_0\end{pmatrix}
 &=\begin{pmatrix}
      y_1\\
      y_2\cosh\tau+y_0\sinh\tau\\
      y_0\cosh\tau+y_2\sinh\tau
    \end{pmatrix}.
 \end{aligned}
\]
The boost changes $u$ to $u+\tau$ and fixes the two sides
of $y_1=\sinh d$. Since
$\langle y'(u),y'(u)\rangle_L=\cosh^2d$, the group
$G_t=B_{t/\cosh d}$ has unit speed on $C$; $d=0$ includes
the geodesic case. In the upper-half-plane metric
$(dx^2+dy^2)/y^2$, a horocycle can be placed at $y=b$.
Translations $G_t(x,y)=(x+bt,y)$ preserve its sides and have
speed $b/y=1$ on that curve.
Conjugation by an ambient isometry gives the required
group for the original $C$.
Choose an orientation along the maximal edge and reinsert
$D_t=G_tD$. The slit $\gamma_t=G_t\gamma$ initially lies in its
relative interior. Continue while
\[
 \operatorname{int}(D_t)\cap N'=\gamma_t^{\circ},\qquad
 \partial D_t\cap N'=\partial\gamma_t,
\]
and the two pieces $G_tD_j,G_tD_k$ lie in the corresponding fixed
$j,k$ sides. Small neighborhoods of the two original endpoints
contain only the smoothed $j/k$ edge of $N'$. Away from those
neighborhoods the compact moved sides are separated from the
fixed network by a positive distance. The preserved intersection
data with $C$ and continuity of $G_t$ therefore make this
admissible condition open initially.
A first new contact occurs at a finite time $t_*$, at latest when
the forward endpoint of $\gamma_t$ reaches the chosen vertex
endpoint of the fixed edge. All moved sets for $0\le t\le t_*$
lie in one common compact region.

Deletion changes the labeled areas by
$(-|D|,|D_j|,|D_k|)$. Isometric reinsertion changes them by the
opposite vector, since
\[
 (|G_tD|,|G_tD_j|,|G_tD_k|)=(|D|,|D_j|,|D_k|).
\]
Every labeled area is therefore preserved before contact.
The interface length is also preserved:
\begin{equation}\label{eq:rolling-length}
 P_t=\mathcal H^1(N')-\mathcal H^1(\gamma_t)
                    +\mathcal H^1(\partial D_t)=P_0.
\end{equation}
For the fixed deleted partition, denote its cells by $E'_h$.
The moving indicators satisfy, up to null sets,
\[
 \mathbf1_{E_i(t)}=\mathbf1_{E'_i}+\mathbf1_{D_t},\quad
 \mathbf1_{E_j(t)}=\mathbf1_{E'_j}-\mathbf1_{G_tD_j},\quad
 \mathbf1_{E_k(t)}=\mathbf1_{E'_k}-\mathbf1_{G_tD_k},
\]
with all other cells fixed. Off the boundaries of the limiting
domains, membership is eventually constant because
$G_t^{-1}x\to G_{t_*}^{-1}x$. These boundaries have area zero,
and the changes are supported in the common compact region.
Dominated convergence gives
\[
 \sum_h\int_{\mathbb H^2}
   |\mathbf1_{E_h(t)}-\mathbf1_{E_h(*)}|\,dA\longrightarrow0.
\]
Their limiting partition $\Lambda_*$ therefore has the prescribed
areas, including the finite area differences for the exterior.

Define its perimeter measure, and the preceding measures, by
\[
 \mu_t=\frac12\sum_h|D\mathbf1_{E_h(t)}|,
 \qquad \mu_*=\frac12\sum_h|D\mathbf1_{E_h(*)}|.
\]
For $t<t_*$ the exact network decomposition gives
\begin{equation}\label{eq:rolling-measures}
 \mu_t=\mathcal H^1\!\llcorner N'
       -\mathcal H^1\!\llcorner\gamma_t
       +\mathcal H^1\!\llcorner\partial D_t.
\end{equation}
For $f\in C_c(\mathbb H^2)$, the area formula for the isometry gives
\[
 \int_{\partial D_t}f\,d\mathcal H^1
 =\int_{\partial D}f(G_tx)\,d\mathcal H^1(x),\qquad
 \int_{\gamma_t}f\,d\mathcal H^1
 =\int_\gamma f(G_tx)\,d\mathcal H^1(x).
\]
The integrands converge uniformly on these fixed compact arcs,
which have finite length. Thus $\mu_t$ converges weakly to the
nonnegative measure
\[
 \widehat\mu_*=\mathcal H^1\!\llcorner N'
       -\mathcal H^1\!\llcorner\gamma_*
       +\mathcal H^1\!\llcorner\partial D_*,\qquad
 \widehat\mu_*(\mathbb H^2)=P_0.
\]
Compact support of all terms gives the total-mass identity.
Weighted lower semicontinuity of each cell's total variation gives,
for $f\in C_c(\mathbb H^2)$, $f\ge0$,
\[
 \begin{split}
 \int f\,d\mu_*
 &=\frac12\sum_h\int f\,d|D\mathbf1_{E_h(*)}|\\
 &\le\frac12\sum_h\liminf_{t\uparrow t_*}
                   \int f\,d|D\mathbf1_{E_h(t)}|\\
 &\le\liminf_{t\uparrow t_*}\int f\,d\mu_t
  =\int f\,d\widehat\mu_*.
 \end{split}
\]
Consequently $\mu_*\le\widehat\mu_*$ as Radon measures.
Fixed-area minimality gives
$\mu_*(\mathbb H^2)\ge P_0$, while lower semicontinuity gives
the reverse inequality. Hence
\begin{equation}\label{eq:rolling-measure-equality}
 P(\Lambda_*)=P_0,\qquad \mu_*=\widehat\mu_*.
\end{equation}
In particular, the limiting partition is again a fixed-area
minimizer and has the planar regularity established in Section~2.

The two moved sides belong to generalized circles different from
$C$. Each such complete circle meets $C$ precisely at its two
prescribed digon endpoints: two distinct circles through two given
points have no third intersection, as follows also by sending
those points to zero and infinity. Positive-length coincidence of
a moved side and a fixed edge must therefore occur off $C$.
On that interval $\widehat\mu_*$ counts both arcs, while the
perimeter measure of the regular limiting partition has
multiplicity at most one. This contradicts
\eqref{eq:rolling-measure-equality}.

For isolated contacts the following two cases exhaust the possibilities.
If $x_*\notin C$, the moving point is in the interior of a moved
side, and contributes two outgoing germs. The fixed film contributes
two germs at a smooth point, or three at a triple vertex. These
arcs are outside the subtracted slit, so at least four germs are
counted by $\widehat\mu_*$.

If $x_*\in C$, a moved side can meet $C$ only at one of its
prescribed endpoints. A new contact there occurs when that
endpoint reaches a fixed vertex of the containing edge. The
$j/k$ segment between the moving and fixed vertices shrinks to
zero. The two non-$C$ germs from the moved vertex and the two
non-$C$ germs from the fixed vertex remain. They again give four
germs in $\widehat\mu_*$. A fixed vertex in the interior of the
moving slit would have been reached by its leading endpoint at
an earlier time. Simultaneous contacts are covered by selecting
any one of their contact points.

For a unit-speed outgoing germ $\gamma(s)$, normal coordinates
at $x_*$ give $\exp_{x_*}^{-1}\gamma(s)=sv+o(s)$ with $|v|=1$.
Its first-exit arclength from $B_r(x_*)$ is therefore
$r+o(r)$. Each of the $N\ge4$ germs contributes at least that
length. Germs with a common tangent still
have additive lengths unless they agree on an interval, which
was excluded above. Consequently
\eqref{eq:rolling-measure-equality} gives
\begin{equation}\label{eq:rolling-contact-density}
 \liminf_{r\downarrow0}\frac{\mu_*(B_r(x_*))}{2r}
 \ge\liminf_{r\downarrow0}\frac{Nr+o(r)}{2r}
 =\frac N2\ge2.
\end{equation}
A regular planar minimizing film has density $1$ at a smooth
point and $3/2$ at a triple point. This contradiction excludes
the first contact and hence the original digon.
\end{proof}

\begin{lemma}[Triangular-face deletion]\label{lem:triangle-deletion}
Let $D$ be a bounded triangular disk face of a regular equilibrium
network. The three external incident arcs, continued inward,
meet at a regular triple point inside $D$ and divide $D$ into
three nonempty open pieces assigned to its former neighbors.
Replacing the old face by these pieces gives an embedded
equilibrium network with the remaining pressures unchanged.
When one neighboring label is $0$, total finite area strictly
decreases. The same area conclusion holds for deletion of a
boundary digon.
\end{lemma}
\begin{proof}
Let the central label be $i$ and its consecutive neighbors be
$j,k,\ell$. These neighbors are pairwise distinct: each pair of
consecutive sides meets at a vertex with three distinct labels.
Work in a Poincar\'e chart, $g_H=e^{2\varphi}g_E$, with the same
normal orientation in both metrics. Let $T_E,N_E$ be Euclidean
unit tangent and normal and set $T_H=e^{-\varphi}T_E$,
$N_H=e^{-\varphi}N_E$. Differentiating the factor $e^{2\varphi}$
in the Koszul formula gives, for arbitrary vector fields $U,V,W$,
\[
 \begin{split}
 2g_H(\nabla^H_UV,W)
 =2e^{2\varphi}\bigl(&\langle D_UV,W\rangle_E
 +(U\varphi)\langle V,W\rangle_E\\
 &+(V\varphi)\langle U,W\rangle_E
 -(W\varphi)\langle U,V\rangle_E\bigr).
 \end{split}
\]
Dividing by $2e^{2\varphi}$ yields
\[
 \nabla^H_UV=D_UV+(U\varphi)V+(V\varphi)U
                       -\langle U,V\rangle_E\nabla_E\varphi.
\]
Using $D_{T_E}N_E=\kappa_ET_E$ and
$\langle T_E,N_E\rangle_E=0$ gives
\[
 \begin{split}
 \nabla^H_{T_H}N_H
 &=e^{-2\varphi}
       \bigl(\nabla^H_{T_E}N_E-(T_E\varphi)N_E\bigr)\\
 &=e^{-2\varphi}(\kappa_E+N_E\varphi)T_E
  =e^{-\varphi}(\kappa_E+N_E\varphi)T_H.
 \end{split}
\]
The convention $\nabla^H_{T_H}N_H=\kappa_HT_H$ therefore gives
\begin{equation}\label{eq:conformal-curve-curvature}
 \kappa_H=e^{-\varphi}
              (\kappa_E+\partial_{N_E}\varphi).
\end{equation}
Each constant-hyperbolic-curvature support is a Euclidean
generalized circle; its signed Euclidean curvature $K^E_{ab}$
is therefore constant along that support.
At each original vertex the three cyclic Euclidean normals sum
to zero. The hyperbolic pressure cocycle and
\eqref{eq:conformal-curve-curvature} give
\[
 \begin{split}
 K^E_{ij}+K^E_{jk}+K^E_{ki}&=0,\\
 K^E_{ik}+K^E_{k\ell}+K^E_{\ell i}&=0,\\
 K^E_{i\ell}+K^E_{\ell j}+K^E_{ji}&=0.
 \end{split}
\]
Adding these equalities yields
\begin{equation}\label{eq:triangle-compatibility}
 K^E_{jk}+K^E_{k\ell}+K^E_{\ell j}=0.
\end{equation}
More explicitly, set $p_i^E=0$ and $p_a^E=-K^E_{ia}$ for
$a=j,k,\ell$. All six local supports then have curvature
$p_a^E-p_b^E$. Thus the central triangular face and its three
incident germs have compatible Euclidean equilibrium data.

Wichiramala's local circular-arc reduction
\cite[Lemmas~5.8--5.9]{Wichiramala} uses the Jordan triangular face,
its $120$-degree angles, and the curvature cocycle at its three
vertices.
Order the vertices positively along $\partial D$ and send them, in
that order, to positively ordered equilateral vertices by a
M\"obius transformation. The circular-arc uniqueness in Lemma~5.8
fixes the selected face and excludes the complementary crossing-arc
choice. The three incident arcs, continued inward, then meet in that
face at one triple point with $120$-degree sectors, as in Lemma~5.9.
Undoing the transformation gives three embedded replacement segments
inside $\overline D$, each on the supporting circle of its old
incident arc. They split $D$ into three nonempty open pieces.

The segments lie in the original compact hyperbolic face and
thus in the Poincar\'e domain. Their unchanged oriented supports
retain the original signed hyperbolic curvatures
$p_j-p_k$, $p_k-p_\ell$, $p_\ell-p_j$.
Angles are conformally preserved. All changes lie inside the
old face, so embeddedness and the exterior network are
preserved. If one recipient is label $0$, its new piece has
positive area and lies inside the old finite union. Total
finite area decreases by exactly the area of that piece.
The digon construction of Lemma~\ref{lem:digon-deletion}
has the same property with two recipient pieces.
\end{proof}

\section{A comparison minimum under a perimeter bound}\label{sec:comparison}

For a hyperbolic disk of radius $R$,
\[
 a=2\pi(\cosh R-1),\qquad I(a)=2\pi\sinh R.
\]
Eliminating $R$ gives
\begin{equation}\label{eq:disk-excess}
 I(a)=\sqrt{a^2+4\pi a},\qquad
 j(a)=I(a)-a,\qquad
 j'(a)=\frac{a+2\pi}{\sqrt{a^2+4\pi a}}-1>0.
\end{equation}
Differentiating once more gives
\begin{equation}\label{eq:disk-excess-concavity}
 j''(a)=\frac{a^2+4\pi a-(a+2\pi)^2}
                    {(a^2+4\pi a)^{3/2}}
        =-\frac{4\pi^2}{(a^2+4\pi a)^{3/2}}<0.
\end{equation}
The positivity of $j'$ follows from
$(a+2\pi)^2-(a^2+4\pi a)=4\pi^2$.
Indeed, $\cosh R=1+a/(2\pi)$ and
\[
 (2\pi\sinh R)^2
 =4\pi^2\bigl((1+a/(2\pi))^2-1\bigr)=a^2+4\pi a.
\]
Lemma~\ref{lem:hyperbolic-isoperimetric} applied to the finite union of
any cluster gives
\begin{equation}\label{eq:positive-excess}
 P(\mathcal C)\ge\operatorname{Per}(E_1\cup E_2\cup E_3)
 \ge I(A)>A,
 \qquad A=v_1+v_2+v_3.
\end{equation}
In particular $J=P-A>0$ and $Q=J/S$ is well defined.

\begin{lemma}[Excluding a vanishing cell]\label{lem:vanishing}
If $v=(u,w,a)>0$ and $Q(\mathcal C)\le\lambda<1$, then
\begin{equation}\label{eq:vanishing}
 (1-\lambda)j(u+w)\le I(a).
\end{equation}
Thus $u+w\ge\delta>0$ implies
\[
 a\ge\sqrt{4\pi^2+(1-\lambda)^2j(\delta)^2}-2\pi>0.
\]
\end{lemma}

\begin{proof}
Relabel the entire third cell as exterior, obtaining a
finite-perimeter double cluster $\Gamma$ of areas $(u,w)$.
Take a perimeter minimizer $\mathcal D$ at these two areas.
The existence and regularity proof of Proposition~\ref{prop:inputs}
applies with two finite labels: the area-correction tree has two
edges, small-cycle erasure saves at least half the cycle length,
and the same local junction argument gives regularity.  Thus
$\mathcal D$ is a bounded circular-arc film and its film length
equals its finite-perimeter cluster perimeter.  Conversely, every
finite-length enclosure in the sense of Cotton--Freeman defines
a finite-perimeter partition whose cluster perimeter is no larger
than its film length.  Hence $\mathcal D$ also minimizes in that
enclosure class.  Their theorem
\cite[Theorem 2.7]{CottonFreeman} identifies it as a standard
double bubble.  Its three circular arcs meet at $120$ degrees,
as do those of our standard two-cell model; the uniqueness in
\cite[Proposition 2.6]{CottonFreeman} identifies the two models.
Consequently
\[
 P(\mathcal C)\ge P(\Gamma)\ge P(\mathcal D)=T_2(u,w).
\]
Writing $S_2=T_2-u-w$ yields
\[
 J\ge S_2-a,\qquad S_3(u,w,a)\le S_2+j(a).
\]
If $S_2-a\ge0$, division gives
$Q\ge(S_2-a)/(S_2+j(a))$; if $S_2-a<0$, the same inequality
follows from $Q>0$.  Therefore
\[
 (1-\lambda)S_2\le a+\lambda j(a)\le I(a).
\]
The disk inequality for the finite union of the standard double
bubble gives $S_2\ge j(u+w)$.  Solving $b^2=a^2+4\pi a$
by completing the square gives $(a+2\pi)^2=4\pi^2+b^2$, and hence
$I^{-1}(b)=\sqrt{4\pi^2+b^2}-2\pi$.
\end{proof}

\begin{lemma}[Uniform small-area comparison]\label{lem:small-area}
Set
\[
 B(A)=\sqrt{3A^2+12\pi A},\qquad
 r(A)=\operatorname{arcosh}\left(1+\frac{A}{2\pi}\right),
\]
\[
 L_m(A)=\cosh^2\frac{B(A)}2,\qquad
 L_s(A)=\cosh^2\frac{3r(A)}2.
\]
Every positive-area three-cluster satisfies
\begin{equation}\label{eq:small-area}
 Q(\mathcal C)\ge
 1-\frac{I(A)}{j(A)}
       \left(1-\frac1{L_m(A)L_s(A)}\right).
\end{equation}
The right-hand side tends to $1$ as $A\downarrow0$.
\end{lemma}

\begin{proof}
Replace $\mathcal C$ by a fixed-area perimeter minimizer
$\mathcal M$.  This decreases $Q$.  The coordinate-face
estimate for the standard family and Cauchy--Schwarz imply
\[
 P(\mathcal M)\le T(v)\le\sum_{i=1}^3I(v_i)
 \le\sqrt{3\sum_{i=1}^3(v_i^2+4\pi v_i)}\le B(A).
\]
The boundary of $\mathcal M$ is a connected finite network.
A simple graph path from a fixed boundary point to any other
boundary point has length at most $P(\mathcal M)$.
Consequently the whole boundary lies in a hyperbolic ball of radius
$B(A)$.  Its bounded faces lie in the same ball: the complement
of that ball is connected and belongs to the unbounded face.

In Poincar\'e coordinates centered at the center of such a ball,
normalized so that the metric agrees with the Euclidean metric at
the origin,
\[
 g_H=\rho^2 g_E,\qquad
 \rho(z)=\frac1{1-|z|^2/4}=\cosh^2\frac{d_H(0,z)}2.
\]
On a ball of radius $R$, $1\le\rho\le\cosh^2(R/2)$.
Hence the Euclidean cell areas $v_i^E$ and Euclidean perimeter
of $\mathcal M$ satisfy
\[
 v_i^E\ge v_i/L_m^2,\qquad P_E(\mathcal M)\le P(\mathcal M).
\]
Let $T_E$ denote the Euclidean standard triple-bubble profile.
The planar triple-bubble theorem \cite{Wichiramala,Lawlor},
coordinatewise monotonicity of $T_E$, and its homogeneity give
\begin{equation}\label{eq:small-actual}
 P(\mathcal M)\ge T_E(v^E)\ge L_m^{-1}T_E(v).
\end{equation}
Indeed, the Euclidean standard cluster
with areas $w$ has exterior defining radii $R_i(w)>0$.
Lawlor's first-variation formula \cite[Proposition 2.13]{Lawlor}
gives $\partial_iT_E(w)=1/R_i(w)>0$ on the positive orthant.
Integration along coordinate segments gives coordinatewise
monotonicity.  Euclidean dilation gives
$T_E(tv)=\sqrt t\,T_E(v)$.

For the hyperbolic standard cluster, each defining disk has area
at most $A$, hence radius at most $r(A)$, by the defining-disk
containment proved in Section~\ref{sec:standard}.  Any two defining
disks intersect, since their cells have a nonempty common interface.
Their union is therefore contained in a ball of radius $3r(A)$
centered at the first disk center.  In its Poincar\'e chart the
standard cluster is an embedded circular-arc $K_4$ equilibrium
network.  Indeed, the conformal curvature formula cancels its
normal-derivative terms at every triple junction.  The four
curvature cocycles on the dual tetrahedron integrate to Euclidean
pressures.  Lemma~\ref{lem:k4} identifies the resulting network
as a Euclidean standard triple bubble.  Its Euclidean areas are
at most $v_i$, whence
\begin{equation}\label{eq:small-model}
 T(v)\le L_s T_E(v^E)\le L_s T_E(v).
\end{equation}
Combining \eqref{eq:small-actual} and \eqref{eq:small-model}
gives $P(\mathcal M)\ge T(v)/(L_mL_s)$.
Since $T(v)\ge I(A)>A$ and $t/(t-A)$ decreases for $t>A$,
\[
 \frac{T(v)}{T(v)-A}\le\frac{I(A)}{I(A)-A}.
\]
Subtracting $A$ and dividing by $S=T-A$ proves
\eqref{eq:small-area}.  Finally $L_m(A),L_s(A)\to1$ and
$I(A)/j(A)\to1$.  These factors depend only on $A$, so the
conclusion is uniform in all three area ratios.
More explicitly, $\cosh^2z=1+z^2+O(z^4)$ gives
\[
 \begin{aligned}
 B(A)^2&=12\pi A+3A^2,
 &L_m(A)&=1+3\pi A+O(A^2),\\
 r(A)^2&=A/\pi+O(A^2),
 &L_s(A)&=1+\frac9{4\pi}A+O(A^2).
 \end{aligned}
\]
Here the formula for $r(A)^2$ follows by substituting
$\cosh r=1+r^2/2+O(r^4)$ into $\cosh r=1+A/(2\pi)$.
Since $I(A)^2-A^2=4\pi A$,
\[
 \frac{I(A)}{j(A)}
 =\frac{I(A)(I(A)+A)}{4\pi A}
 =1+\frac{A+I(A)}{4\pi}=1+O(\sqrt A).
\]
Thus the right-hand side of \eqref{eq:small-area} equals
\[
 1-\left(3\pi+\frac9{4\pi}\right)A+O(A^{3/2}).
\]
\end{proof}

\begin{corollary}\label{cor:positive-coordinates}
For every $\lambda<1$ there is $\varepsilon_\lambda>0$ such
that $Q(\mathcal C)\le\lambda$ implies
$v_i\ge\varepsilon_\lambda$ for $i=1,2,3$.
\end{corollary}

\begin{proof}
Choose $\delta_\lambda>0$ so that the right-hand side of
\eqref{eq:small-area} exceeds $\lambda$ for
$0<A<\delta_\lambda$.  Thus $A\ge\delta_\lambda$.
If $v_i<\delta_\lambda/2$, the other two coordinates sum to
at least $\delta_\lambda/2$ and Lemma~\ref{lem:vanishing}
applies.  One can take
\[
 \varepsilon_\lambda=\min\left\{
 \frac{\delta_\lambda}{2},
 \sqrt{4\pi^2+(1-\lambda)^2j(\delta_\lambda/2)^2}-2\pi
 \right\}.
\]
\end{proof}

\begin{proposition}[Attainment under a perimeter bound]\label{prop:attainment}
If $Q(\mathcal C_0)=\lambda_0<1$, then for every
$M\ge P(\mathcal C_0)$ the ratio $Q$ attains a minimum
$0<\lambda\le\lambda_0$ among positive-area clusters with
$P\le M$.  One such minimizer is a perimeter minimizer for its
own area vector.
\end{proposition}

\begin{proof}
Define
\[
 q(v)=\frac{I_3(v)-A(v)}{S(v)},\qquad
 K=\{v\in(0,\infty)^3:I_3(v)\le M,\ q(v)\le\lambda_0\}.
\]
Replacing $\mathcal C_0$ by a fixed-area minimizer shows that
$K\ne\varnothing$.  Inequality~\eqref{eq:positive-excess}
and Corollary~\ref{cor:positive-coordinates} place $K$ inside
$[\varepsilon_{\lambda_0},M]^3$.  The continuity of $I_3$ and
$T$ makes $K$ closed in this box.  Thus $q$ attains its minimum
on $K$.  Its value is positive by \eqref{eq:positive-excess}.
Choose a fixed-area minimizer at the minimizing vector.  A
cluster under the same perimeter bound with a smaller value of
$Q$ could itself be replaced by a fixed-area minimizer.  Its
area vector would belong to $K$, contradicting the minimum of $q$.
\end{proof}

\begin{lemma}[Pressure relation]\label{lem:capped-pressure}
At a comparison minimum supplied by Proposition~\ref{prop:attainment},
let $\lambda=Q(\mathcal C_*)<1$ and let $\kappa$ be the model
pressure vector at its areas.  There is $\eta\ge0$ such that
\begin{equation}\label{eq:capped-pressure}
 p_i=\frac{1+\lambda(\kappa_i-1)}{1+\eta S},
 \qquad 0<p_i<\kappa_i.
\end{equation}
If $P(\mathcal C_*)<M$, then $\eta=0$.
\end{lemma}

\begin{proof}
For every compactly supported smooth field $X$,
\[
 \begin{split}
 DQ[X]
 &=\frac{(DP[X]-\sum_iDv_i[X])S
       -J\sum_i(\kappa_i-1)Dv_i[X]}{S^2}\\
 &=\frac{DP[X]-\sum_i[1+\lambda(\kappa_i-1)]Dv_i[X]}S.
 \end{split}
\]
When $P<M$, both signs of each sufficiently small variation are
admissible, so $DQ=0$.  When $P=M$, take a compactly supported
field equal to $X_0=-\nabla(r^2/2)$ near the cluster, where
$r$ is distance from a fixed point.  In hyperbolic polar coordinates,
$g=dr^2+\sinh^2r\,d\theta^2$.  The connection coefficient
$\Gamma^r_{\theta\theta}=-\sinh r\cosh r$ gives
\[
 \nabla^2r=\coth r\,(g-dr\otimes dr),\qquad
 \nabla^2\frac{r^2}{2}
 =dr\otimes dr+r\coth r\,(g-dr\otimes dr).
\]
The squared distance is smooth also at the center.  Therefore
\[
 \operatorname{div}_{\tau}X_0
 =-\bigl(dr(\tau)^2+r\coth r\,[1-dr(\tau)^2]\bigr)\le-1
\]
for every unit tangent $\tau$, using the limiting value $1$ at
$r=0$.  Therefore $DP[X_0]\le-P<0$.
Every $X$ with $DP[X]<0$ is admissible for small positive time
and satisfies $DQ[X]\ge0$.  If $DP[Y]=0$, apply this observation
to $Y+\varepsilon X_0$ and $-Y+\varepsilon X_0$, then let
$\varepsilon\downarrow0$.  It follows that $DQ[Y]=0$.
Thus $DQ=-\eta DP$, where
$\eta=-DQ[X_0]/DP[X_0]\ge0$.  Substitution gives
\[
 (1+\eta S)DP=\sum_i[1+\lambda(\kappa_i-1)]Dv_i.
\]
Lemma~\ref{lem:area-surjectivity} identifies the coefficients,
yielding \eqref{eq:capped-pressure}.  Its numerator is positive
and strictly smaller than $\kappa_i$ because
$0<\lambda<1<\kappa_i$; its denominator is at least one.
\end{proof}

\section{Relabeling estimates and planar component restrictions}
\label{sec:relabeling}\label{sec:dimensional-relabeling}

We first establish the relabeling estimates in $\mathbb H^n$, then
specialize to the three-bubble problem in $\mathbb H^2$.
Fix $n\ge2$, put $m=n-1$ and $\alpha=m/n$, and retain the ball
profile $B_n$ and standard excess $S_{n,k}$ from
Section~\ref{sec:standard}.

\begin{lemma}[Single-ball excess]\label{lem:single-ball-excess}
For $a>0$, the function $j_n(a)=B_n(a)-ma$ satisfies
$j_n(a)>0$, $j_n'(a)>0$, and $j_n(a)/a^\alpha$ is strictly decreasing.
Consequently, if $0<b<a$, then
\begin{equation}\label{eq:dimensional-excess-loss}
 j_n(a)-j_n(a-b)<\Phi_\alpha(b/a)j_n(b),\qquad
 \Phi_\alpha(r)=\frac{1-(1-r)^\alpha}{r^\alpha}.
\end{equation}
\end{lemma}
\begin{proof}
Let $\omega=|\mathbb S^{n-1}|$ and write
$a=\omega\int_0^R\sinh^m t\,dt$. Then
\[
 B_n(a)=\omega\sinh^mR,\qquad
 j_n'(a)=m(\coth R-1)>0.
\]
Since $j_n(0)=0$, positivity follows. Since $\alpha=m/n$,
\[
 \begin{split}
 \frac{d}{da}\left(\frac{j_n(a)}{a^\alpha}\right)
 &=a^{-\alpha-1}\bigl(a j_n'(a)-\alpha j_n(a)\bigr),\\
 a j_n'(a)-\alpha j_n(a)
 &=ma(\coth R-1)-\frac mn(B_n(a)-ma)\\
 &=\frac mn\bigl(a(n\coth R-1)-B_n(a)\bigr).
 \end{split}
\]
Thus the derivative is negative precisely when
$B_n(a)>a(n\coth R-1)$, equivalently
$j_n(a)>na(\coth R-1)$. Substitution reduces this
inequality to
\[
 \int_0^R\sinh^m t\,dt<h(R),\qquad
 h(R)=\frac{\sinh^mR}{n\coth R-1}.
\]
The denominator is positive. As $R\downarrow0$,
$\sinh R=R+O(R^3)$ and $\coth R=R^{-1}+O(R)$, so
$h(R)=R^{m+1}/n+O(R^{m+2})\to0$. Put $c=\coth R$ and
$U=\sinh^mR$, so $c'=1-c^2$ and $U'=mcU$. Then
\[
 h'(R)=\frac{mcU(nc-1)+nU(c^2-1)}{(nc-1)^2},
\]
and therefore
\[
 \begin{split}
 h'(R)-U
 &=\frac{U}{(nc-1)^2}
       \bigl(mc(nc-1)+n(c^2-1)-(nc-1)^2\bigr)\\
 &=\frac{U}{(nc-1)^2}
       \bigl((mn+n-n^2)c^2+(2n-m)c-(n+1)\bigr)\\
 &=\frac{(n+1)(c-1)U}{(nc-1)^2},
 \end{split}
\]
where $m=n-1$ was used in the last line. Hence
\[
 h'(R)-\sinh^mR
  =\frac{(n+1)(\coth R-1)\sinh^mR}{(n\coth R-1)^2}>0.
\]
Integration from zero proves the required strict inequality.
The strictly decreasing quotient gives
\[
 j_n(a-b)>\Bigl(\frac{a-b}{a}\Bigr)^\alpha j_n(a),\qquad
 j_n(b)>\Bigl(\frac ba\Bigr)^\alpha j_n(a).
\]
Multiplying the second inequality by $\Phi_\alpha(b/a)>0$ and
adding the first proves \eqref{eq:dimensional-excess-loss}.
\end{proof}

For $2\le k\le n+1$, a finite-perimeter $k$-cluster is a partition
$\mathcal C=(E_0,\ldots,E_k)$ of $\mathbb H^n$ into cells of locally
finite perimeter, with $0<v_i=|E_i|<\infty$ for $i\ge1$ and
$P(\mathcal C)=\frac12\sum_{i=0}^k\operatorname{Per}(E_i)<\infty$.
Define
\[
 J_n(\mathcal C)=P(\mathcal C)-m\sum_i v_i,\qquad
 Q_n(\mathcal C)=\frac{J_n(\mathcal C)}{S_{n,k}(v)}.
\]
The denominator is positive by Corollary~\ref{cor:standard-excess},
and Lemma~\ref{lem:hyperbolic-isoperimetric} applied to the finite
union gives $J_n\ge j_n(\sum_i v_i)>0$.

\begin{lemma}[Relabeling identities]
\label{lem:dimensional-relabeling-accounting}
Let $K\subset E_i$ up to a null set be a measurable set of finite
perimeter and volume $b=|K|$, where $0<b<a=v_i$. Suppose
\begin{equation}\label{eq:interface-exhaustion}
 \ell_j=\mathcal H^{n-1}(\partial^*K\cap\partial^*E_j)\quad(j\ne i),
 \qquad \operatorname{Per}(K)=\sum_{j\ne i}\ell_j.
\end{equation}
Relabel $K$ by any $s\ne i$, and let $\mathcal C'$ and $v'$ denote
the resulting cluster and volume vector. Then
$P(\mathcal C')=P(\mathcal C)-\ell_s$, every finite volume remains
positive, and the decrease $d_s=J_n(\mathcal C)-J_n(\mathcal C')$ is
\[
 d_s=\ell_s\quad(s\ne0,i),\qquad d_0=\ell_0-mb.
\]
Putting $D=S_{n,k}(v)-S_{n,k}(v')$, one has
\begin{equation}\label{eq:dimensional-loss-bound}
 D\le j_n(a)-j_n(a-b)<\Phi_\alpha(b/a)j_n(b),
\end{equation}
and
\begin{equation}\label{eq:dimensional-relabeling-accounting}
 Q_n(\mathcal C')-Q_n(\mathcal C)
   =\frac{Q_n(\mathcal C)D-d_s}{S_{n,k}(v')}.
\end{equation}
\end{lemma}
\begin{proof}
Set $R=E_i\setminus K$ and refine the partition by replacing $E_i$
with $K,R$. Both are of finite perimeter. The set-operation and
cluster-interface identities of
\cite[Theorem 16.3 and Proposition 29.4]{Maggi} give, up to
$\mathcal H^{n-1}$-null sets,
\[
 \begin{gathered}
 \operatorname{Per}(K)
 =\mathcal H^{n-1}(\partial^*K\cap\partial^*R)
     +\sum_{j\ne i}\ell_j,\\
 \mathcal H^{n-1}(\partial^*K\cap\partial^*R)=0,\\
 P(\widehat{\mathcal C})=P(\mathcal C),\qquad
 P(\mathcal C')=P(\widehat{\mathcal C})-\ell_s
              =P(\mathcal C)-\ell_s.
 \end{gathered}
\]
The second line uses \eqref{eq:interface-exhaustion}; the last merges
$K$ into the recipient cell in the refined partition
$\widehat{\mathcal C}$. These local BV identities hold in smooth
coordinate charts and hence for hyperbolic perimeter. In particular,
no interface with $R$ is introduced. A finite recipient gains volume
$b$, while recipient $0$ removes $b$ from total finite volume.
This proves the formulas for $d_s$. Since $b<a$, all finite
volumes stay positive, so $S_{n,k}(v')>0$.

First decrease coordinate $i$ by $b$. Corollary~\ref{cor:standard-excess}
gives
\[
 S_{n,k}(v)-S_{n,k}(v-be_i)
 \le\int_{a-b}^{a}j_n'(t)\,dt=j_n(a)-j_n(a-b).
\]
Increasing a finite recipient coordinate afterwards increases
$S_{n,k}$; for recipient $0$ there is no additional change.
Lemma~\ref{lem:single-ball-excess} proves
\eqref{eq:dimensional-loss-bound}. Direct subtraction gives
\eqref{eq:dimensional-relabeling-accounting}.
\end{proof}

\begin{proposition}[Component relabeling]
\label{prop:dimensional-relabeling}
In the setting of Lemma~\ref{lem:dimensional-relabeling-accounting},
suppose $0<Q_n(\mathcal C)\le1$.
If $\ell_0=0$ and $\Phi_\alpha(b/a)\le1/(k-1)$, relabeling $K$ to
some other finite label strictly decreases $Q_n$. If $\ell_0>0$ and
$\Phi_\alpha(b/a)\le1/k$, relabeling it to some other label, possibly
the exterior, strictly decreases $Q_n$. Each relabeling preserves
positivity of the finite-cell volumes and does not increase perimeter.
\end{proposition}
\begin{proof}
If $\ell_0=0$, the $k-1$ finite-recipient savings sum to
$\operatorname{Per}(K)$. One of them satisfies
\[
 d_s\ge\frac{\operatorname{Per}(K)}{k-1}
 \ge\frac{B_n(b)}{k-1}>\frac{j_n(b)}{k-1}.
\]
If $\ell_0>0$, sum all $k$ possible recipient savings, including $0$:
\[
 \sum_{s\ne i}d_s=\operatorname{Per}(K)-mb\ge j_n(b).
\]
Some operation thus has $d_s\ge j_n(b)/k>0$. Under the corresponding
threshold hypothesis, \eqref{eq:dimensional-loss-bound} gives
$D<d_s$ in either case. Put $q=Q_n(\mathcal C)$.
If $D\ge0$, then $qD\le D<d_s$ because $q\le1$; if $D<0$,
then $qD<0<d_s$ because $q>0$.
Equation~\eqref{eq:dimensional-relabeling-accounting} therefore
gives a strict decrease in both cases. The remaining assertions
follow from Lemma~\ref{lem:dimensional-relabeling-accounting}.
\end{proof}
For a regular planar cluster, every whole component of an open finite
cell satisfies \eqref{eq:interface-exhaustion}.

The threshold function is strictly increasing, since
\[
 \begin{split}
 \Phi_\alpha'(r)
 &=-\alpha r^{-\alpha-1}\bigl(1-(1-r)^\alpha\bigr)
       +\alpha r^{-\alpha}(1-r)^{\alpha-1}\\
 &=\alpha r^{-\alpha-1}\bigl((1-r)^{\alpha-1}-1\bigr)>0.
 \end{split}
\]
The expansion $1-(1-r)^\alpha=\alpha r+O(r^2)$ gives
$\Phi_\alpha(0+)=0$, while $\Phi_\alpha(1-)=1$.
Thus every level in $(0,1)$ has a unique threshold in $(0,1)$.
For $k=2$ the interior level is $1$, so every proper subset
satisfying interface exhaustion and $\ell_0=0$ meets its threshold.
For $n=2$, $\alpha=1/2$. If $0<\gamma<1$, then
\[
 \Phi_{1/2}(r)=\gamma
 \ \Longleftrightarrow\
 \sqrt{1-r}=1-\gamma\sqrt r
 \ \Longrightarrow\
 (1+\gamma^2)r=2\gamma\sqrt r.
\]
Since $r>0$, the threshold is
\[
 r=\frac{4\gamma^2}{(1+\gamma^2)^2}.
\]
Conversely, this value has
$1-\gamma\sqrt r=(1-\gamma^2)/(1+\gamma^2)>0$, so squaring
introduced no extraneous root. The $k=3$ interior and
exterior-interface thresholds are therefore $16/25$ and $9/25$.

We now specialize to $n=2$, $k=3$. Then $m=1$, $\alpha=1/2$,
$B_2=I$, $j_2=j$, $S_{2,3}=S$, $J_2=J$, and $Q_2=Q$.
The next two statements record the planar corollaries needed below.

\begin{lemma}[Planar disk-excess estimate]\label{lem:disk-loss}
For $0<b<a$,
\begin{equation}\label{eq:disk-loss}
 j(a)-j(a-b)<c(b/a)j(b),\qquad
 c(t)=\frac{1-\sqrt{1-t}}{\sqrt t}
     =\frac{\sqrt t}{1+\sqrt{1-t}}.
\end{equation}
The function $c$ increases on $(0,1)$, and
$c(16/25)=1/2$, $c(9/25)=1/3$.
\end{lemma}
\begin{proof}
Take $n=2$ in Lemma~\ref{lem:single-ball-excess}; then
$c=\Phi_{1/2}$. The preceding threshold computation gives the
two stated values.
\end{proof}

\begin{lemma}[Planar relabeling identities]\label{lem:relabel}
Let $\mathcal C$ be a regular cluster with compact boundary and let $K$ be a
whole connected component of finite label $i$. Write
$b=|K|$, $a=v_i$, where $0<b<a$. For $s\ne i$ let
$\ell_s$ be the interface length between $K$ and label $s$.
Relabeling $K$ by a finite label $s$ decreases $J$ by
$d_s=\ell_s$; relabeling it by the exterior decreases $J$ by
$d_0=\ell_0-b$. If $v'$ is the resulting area vector and
$D=S(v)-S(v')$, then
\begin{equation}\label{eq:relabel}
 D\le j(a)-j(a-b),\qquad
 Q'-Q=\frac{QD-d_s}{S(v')}.
\end{equation}
Each relabeling does not increase perimeter, and all finite cells
retain positive area.
\end{lemma}
\begin{proof}
The boundary interfaces of a whole regular component exhaust its
perimeter. Apply Lemma~\ref{lem:dimensional-relabeling-accounting}
with $n=2$, $k=3$, and $m=1$.
\end{proof}

\begin{proposition}[The $36\%/64\%$ restriction]
\label{prop:components}
Suppose $C$ is a regular cluster with compact boundary minimizing $Q$ in a class
closed under perimeter-nonincreasing whole-component relabelings,
and $0<Q(C)\le1$.  Its exterior is connected.  Each finite label has
either one interior component and no boundary component, or one or
two boundary components and no interior component.
\end{proposition}
\begin{proof}
If a bounded exterior component $K$ has area $b$, choose an adjacent
finite label $i$ with $\ell_i>0$ and relabel $K$ by $i$.  Then
\[
 P'=P-\ell_i,\quad A'=A+b,\quad J'=J-\ell_i-b<J,
 \qquad S(v+be_i)>S(v).
\]
Both new excesses are positive, so $Q'<Q$, a contradiction.  The
bounded network has exactly one unbounded complementary component;
therefore the exterior has just that component.

Let $K$ be a component of finite label $i$ with
$0<b=|K|<a=v_i$. Its boundary satisfies
\eqref{eq:interface-exhaustion}. Apply
Proposition~\ref{prop:dimensional-relabeling} with $n=2$, $k=3$.
By Lemma~\ref{lem:disk-loss}, $c=\Phi_{1/2}$ is increasing, with $c(16/25)=1/2$ and
$c(9/25)=1/3$. Consequently, an interior component with $b/a\le16/25$ or a
boundary component with $b/a\le9/25$ admits a strict decrease of
$Q$ by an admissible relabeling. Minimality excludes both cases.
Thus, whenever a label has more than one component, every interior
component has area greater than $16/25$ of the label's total area and every
boundary component has area greater than $9/25$ of the total.
Two interior components cannot coexist.  An interior and a boundary
component cannot coexist, since their proportions would sum to more
than one.  Three boundary components cannot coexist either.  A label
with one component already has the asserted form.
\end{proof}

\begin{proposition}[The two combinatorial alternatives]
\label{prop:two-types}
Suppose a cluster has the compact connected regular network of
Lemma~\ref{lem:disk-faces}, connected exterior, no bounded digon,
and the component restrictions of Proposition~\ref{prop:components}.
Either every finite component meets the exterior, or the network is
a 4-daisy: a single interior component surrounded by four
alternating boundary components of the other two labels.
\end{proposition}
\begin{proof}
At least one finite label meets the exterior.  If the other two
labels had only interior components, there could be no exterior
triple junction: such a junction has two distinct finite labels.
The outer boundary would then be a complete vertex-free circle and,
by connectedness of the network, would be the whole film.  This
cannot enclose three finite cells of positive area.  Hence at most one
finite cell has an interior component.

Suppose label $3$ has its unique interior component $D$.  The other
two labels have at most four components in total, each meeting the
exterior.  Along $\partial D$ their labels $1,2$ alternate at every
vertex.  The number $d$ of sides of $D$ is even and, since there is
no digon, $d\ge4$.

No boundary component $E$ can share two distinct edges with $D$.
Indeed choose points $P,Q$ in the interiors of two such edges.
Because $D$ and $E$ are Jordan disks, join $P$ to $Q$ by a proper
crosscut through each disk.  Their union is a Jordan curve $\Gamma$,
otherwise disjoint from the film.  Suppose $E$ has label $1$.
Each of the two portions of $\partial D$ between the chosen
$1/3$ edges contains a $2/3$ edge.  Small label-$2$ patches across
these edges lie on opposite sides of $\Gamma$: the crosscut in $D$
separates its two boundary portions, and the crosscut in $E$ is
disjoint from them.  The curve $\Gamma$ lies in labels $1,3$ and
the chosen interior $1/3$ points, so it misses both label $2$ and
the exterior.  The connected unbounded exterior lies on its
unbounded side.  The label-$2$ component on its bounded side cannot
cross $\Gamma$ and hence cannot meet the exterior.  This contradicts the fact that every component of label $2$
meets the exterior.

Thus each neighbor of $D$ appears once.  If $B$ is the number of
boundary components, then $4\le d\le B\le4$.  Therefore $d=B=4$,
and both boundary labels have two components.  There are six faces, including $D$ and the exterior.  Euler's formula and trivalence
give eight vertices and twelve edges.  Four vertices lie on $D$;
the remaining four lie on the exterior, since every vertex involves
three distinct labels and no other interior face exists.  There are
four exterior edges, and every boundary component has at least one,
so each has exactly one.  The four consecutive pairs of petals around
$D$ are different.  A petal--petal edge starting at a central vertex
therefore cannot end at another central vertex, since its two
faces adjacent to that edge would have to be the same at both ends.
It ends at the exterior.  These four edges, the four central sides
and the four exterior sides exhaust the twelve edges, giving the
4-daisy.
\end{proof}

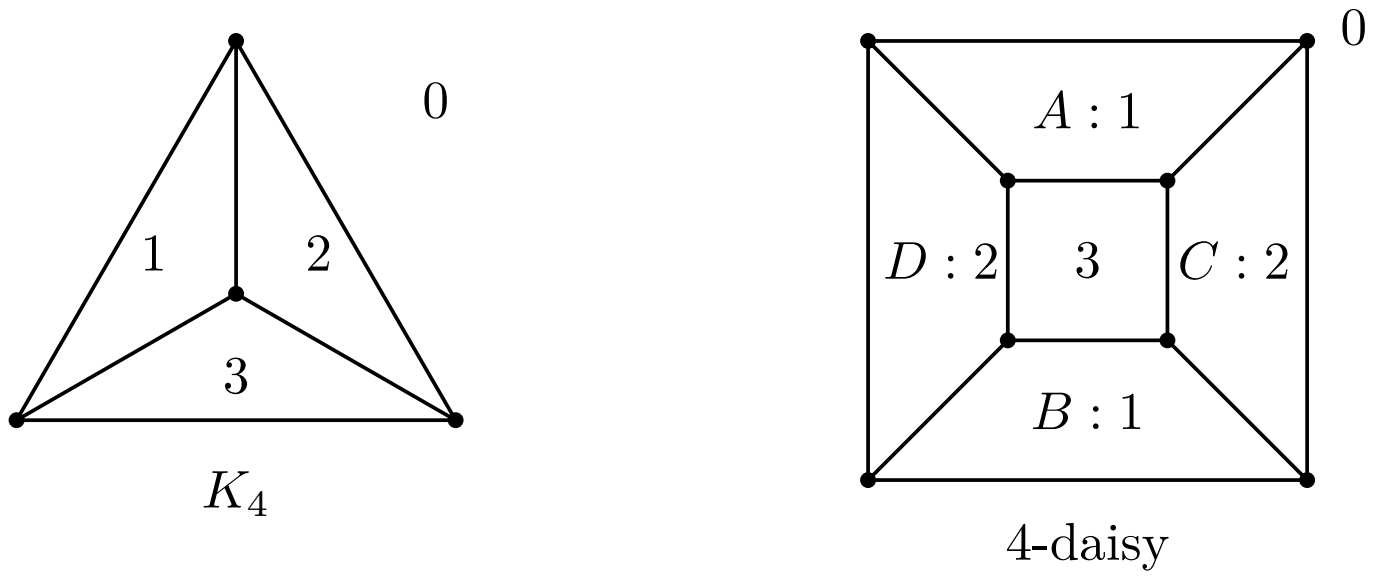
\begin{figure}[t]
\centering
\begin{tikzpicture}[scale=.88,line width=.7pt,every node/.style={font=\small}]
\begin{scope}[xshift=-3.4cm]
\coordinate (a) at (0,1.65);
\coordinate (b) at (-1.65,-1.2);
\coordinate (c) at (1.65,-1.2);
\coordinate (d) at (0,-.25);
\draw (a)--(b)--(c)--cycle;
\draw (a)--(d)--(b) (d)--(c);
\foreach \p in {a,b,c,d}\fill (\p) circle (1.7pt);
\node at (-.62,.05) {$1$};
\node at (.62,.05) {$2$};
\node at (0,-.87) {$3$};
\node at (1.5,1.2) {$0$};
\node at (0,-1.75) {$K_4$};
\end{scope}
\begin{scope}[xshift=3cm]
\coordinate (a) at (-.6,.6);
\coordinate (b) at (.6,.6);
\coordinate (c) at (.6,-.6);
\coordinate (d) at (-.6,-.6);
\coordinate (e) at (-1.65,1.65);
\coordinate (f) at (1.65,1.65);
\coordinate (g) at (1.65,-1.65);
\coordinate (h) at (-1.65,-1.65);
\draw (a)--(b)--(c)--(d)--cycle;
\draw (e)--(f)--(g)--(h)--cycle;
\draw (a)--(e) (b)--(f) (c)--(g) (d)--(h);
\foreach \p in {a,b,c,d,e,f,g,h}\fill (\p) circle (1.7pt);
\node at (0,0) {$3$};
\node at (0,1.13) {$A:1$};
\node at (0,-1.13) {$B:1$};
\node at (1.1,0) {$C:2$};
\node at (-1.1,0) {$D:2$};
\node at (2,1.75) {$0$};
\node at (0,-2.15) {4-daisy};
\end{scope}
\end{tikzpicture}
\caption{The tetrahedral network and the 4-daisy, shown as plane
graphs.  Numerals denote cell labels.}
\label{fig:types}
\end{figure}

\section{Excluding the 4-daisy}

\begin{proposition}[A strict relabeling improvement]\label{prop:no-daisy}
Every 4-daisy with $0<Q\le1$ admits a whole-petal relabeling
which does not increase perimeter, preserves all three finite labels,
and strictly decreases $Q$.
\end{proposition}
\begin{proof}
Write the petals in cyclic order as $A,C,B,D$.  The areas of
label-$1$ petals $A,B$ are $a,b$, those of label-$2$ petals $C,D$
are $c,d$, and the center has label $3$ and area $w>0$.  Write
$\ell_{UV}$ for the interface length between faces $U,V$.
The initial area vector is $(a+b,c+d,w)$.  Renaming the two boundary
labels and their petals, arrange that $a$ is largest and $c\ge d$.

If $b>d$, interchange the labels of $B,D$.  The difference of the
two label totals changes from $(a-c)+(b-d)$ to
$(a-c)-(b-d)$, whose absolute value is no larger.  At a fixed sum
and fixed third coordinate, symmetry and concavity of $S$ imply
that this more balanced pair has no smaller value of $S$: it is a
convex combination of the original pair and its coordinate
interchange.  The switch erases positive-length petal interfaces,
strictly decreasing $P$ and $J$, while preserving total finite area.
More explicitly, put $\delta=(a-c)+(b-d)>0$ and
$\delta'=(a-c)-(b-d)$.  Since $|\delta'|\le\delta$,
$t=(1+\delta'/\delta)/2$ belongs to $[0,1]$, and
\[
 (a+d,b+c)=t(a+b,c+d)+(1-t)(c+d,a+b).
\]
Concavity and symmetry give the claimed inequality for $S$.
The erased interfaces are $AD$ and $BC$, so the actual saving is
$\ell_{AD}+\ell_{BC}>0$.  This decreases $Q$.  We may assume
\[
 a\ge c\ge d\ge b.
\]
If $b\le9a/16$, then $b/(a+b)\le9/25$ and the boundary-component
comparison in Proposition~\ref{prop:components} already improves
$Q$.  Hence assume $b>9a/16$.

Let
\[
 x=\ell_{AC}+\ell_{CB}+\ell_{BD}+\ell_{DA}>0,\qquad
 H=j(a+b)-j(2b)\ge0,
\]
where $2b\le a+b$ and $j$ is increasing.  Each comparison starts from
the original cluster.  Write $d_r=J-J_r$ and $D_r=S-S_r$; the sign
of its ratio change is the sign of $QD_r-d_r$.

Suppose first that $x>H$.  If $a+b\ge c+d$, separately interchange
$B,C$ and $B,D$.  Every petal--petal interface is erased in exactly
one of these two comparisons.  The precise bookkeeping is
\[
\begin{array}{c|c|c|c}
 \text{switch}&v'&d_r&D_r\text{ upper bound}\\\hline
 B,C&(a+c,b+d,w)&\ell_{AC}+\ell_{BD}&j(c+d)-j(b+d)\\
 B,D&(a+d,b+c,w)&\ell_{AD}+\ell_{BC}&j(c+d)-j(b+c)
\end{array}
\]
For the first row, coordinate $2$ loses $c-b$; for the second it
loses $d-b$.  Shrinking that coordinate first and then increasing
coordinate $1$ proves the displayed bounds exactly as in
Lemma~\ref{lem:relabel}.  Thus $d_1+d_2=x$ and
\[
 D_1+D_2\le2j(c+d)-j(b+d)-j(b+c)
             \le j(c+d)-j(2b)\le H.
\]
For the middle inequality, put $s=c+d$.  Since $c,d\ge b$, the
two arguments $b+c,b+d$ have sum $s+2b$ and each belongs to
$[2b,s]$.  Concavity bounds their $j$-sum below by
$j(2b)+j(s)$, its value at the most separated pair.

If $c+d\ge a+b$, separately interchange $A,C$ and $A,D$.
Here the corresponding data are
\[
\begin{array}{c|c|c|c}
 \text{switch}&v'&d_r&D_r\text{ upper bound}\\\hline
 A,C&(b+c,a+d,w)&\ell_{BC}+\ell_{AD}&j(a+b)-j(b+c)\\
 A,D&(b+d,a+c,w)&\ell_{AC}+\ell_{BD}&j(a+b)-j(b+d)
\end{array}
\]
Coordinate $1$ loses $a-c$ or $a-d$, respectively.  It is shrunk
first; the subsequent increase in coordinate $2$ increases $S$.
Again the savings sum to $x$, and
\[
 D_1+D_2\le2j(a+b)-j(b+c)-j(b+d)\le H.
\]
To verify the last inequality, decrease $c,d$ coordinatewise within
$[b,a]$ to $c',d'$ with $c'+d'=a+b$; such a pair exists because
$c+d\ge a+b\ge2b$.  Monotonicity and then concavity give
\[
 j(b+c)+j(b+d)\ge j(b+c')+j(b+d')
                         \ge j(a+b)+j(2b).
\]
In either subcase $\sum d_r=x>H\ge\sum D_r$.  Since $H\ge0$
and $0<Q\le1$, we have $Q\sum D_r\le H$ whether
$\sum D_r$ is nonnegative or negative.  Hence
$\sum(QD_r-d_r)<0$, so one comparison strictly improves the ratio.

Now suppose $x\le H$.  Make four separate comparisons: relabel
$A$ by the center, relabel $A$ by the exterior, and do the same
two operations for $B$.  Their complete data are
\[
\begin{array}{c|c|c|c}
 \text{relabeling}&v'&d_r&D_r\text{ upper bound}\\\hline
 A\to3&(b,c+d,w+a)&\ell_{A3}&j(a+b)-j(b)\\
 A\to0&(b,c+d,w)&\ell_{A0}-a&j(a+b)-j(b)\\
 B\to3&(a,c+d,w+b)&\ell_{B3}&j(a+b)-j(a)\\
 B\to0&(a,c+d,w)&\ell_{B0}-b&j(a+b)-j(a)
\end{array}
\]
Since
\[
 \begin{split}
 \operatorname{Per}(A)&=\ell_{A3}+\ell_{A0}+\ell_{AC}+\ell_{AD},\\
 \operatorname{Per}(B)&=\ell_{B3}+\ell_{B0}+\ell_{BC}+\ell_{BD},
 \end{split}
\]
their actual excess savings satisfy exactly
\[
 \begin{split}
 \sum_r d_r
 &=\operatorname{Per}(A)+\operatorname{Per}(B)-x-a-b\\
 &\ge j(a)+j(b)-x\\
 &\ge j(a)+j(b)-j(a+b)+j(2b).
 \end{split}
\]
Here $-a,-b$ are the area corrections for the two exterior
relabelings.  The model savings satisfy
\[
 \sum_r D_r\le4j(a+b)-2j(a)-2j(b)=:U,
\]
where
\[
 U=2\bigl(j(a+b)-j(a)\bigr)+2\bigl(j(a+b)-j(b)\bigr)\ge0.
\]
We show that the preceding lower bound for $\sum d_r$ is strictly
larger than $U$, equivalently
\begin{equation}\label{eq:daisy-strict}
 3j(a)+3j(b)+j(2b)>5j(a+b).
\end{equation}
Set $\theta=b/(a+b)\in(9/25,1/2]$.  Since $j(t)/\sqrt t$
decreases and $a,b,2b\le a+b$,
\[
 \frac{3j(a)+3j(b)+j(2b)}{j(a+b)}
 \ge f(\theta):=3\sqrt{1-\theta}+(3+\sqrt2)\sqrt\theta.
\]
On $[9/25,1/2]$,
\[
 f''(\theta)=-\frac3{4(1-\theta)^{3/2}}
             -\frac{3+\sqrt2}{4\theta^{3/2}}<0.
\]
Thus $f$ is bounded below by the smaller endpoint value, and both
\[
 \begin{split}
 f(9/25)&=3\frac45+(3+\sqrt2)\frac35
          =\frac{21+3\sqrt2}{5}
          =5+\frac{3\sqrt2-4}{5}>5,\\
 f(1/2)&=\frac{3+3+\sqrt2}{\sqrt2}
          =1+3\sqrt2=5+(3\sqrt2-4)>5.
 \end{split}
\]
Indeed $3\sqrt2>4$ because both sides are positive and
$(3\sqrt2)^2=18>16=4^2$.
This proves \eqref{eq:daisy-strict}, hence $\sum d_r>U$.
Again $Q\sum D_r\le U$, treating a negative sum by
$Q\sum D_r<0\le U$.  Some $QD_r-d_r$ is negative.
Each new denominator is positive, so that comparison improves $Q$.

All operations preserve positive area in each finite cell: a second
petal of each boundary label survives and the center is not erased.
They only change labels on whole petals, erasing interfaces and
creating none.  Thus any prescribed perimeter bound is preserved.
\end{proof}

\section{Pressure-preserving pruning and the conclusion}

\begin{lemma}[Pruning an all-boundary network]\label{lem:pruning}
Let a regular equilibrium cluster with compact boundary have connected boundary and
connected exterior.  Suppose every finite face meets the exterior,
there is no bounded digon, and all three finite cells have positive
area.  A finite sequence of triangular-face and digon deletions
preserves all three labels and their pressures and ends at a
$K_4$ network.  If a deletion occurs, the endpoint has strictly smaller
total finite area.
\end{lemma}
\begin{proof}
Call the original embedded network the primal network.
The elementary graph argument of Lemma~\ref{lem:disk-faces} applies
to this regular network: an actual interface separates distinct
labels, so no bridge occurs; trivalence then excludes loops and cut
vertices.  All primal faces on the compactified plane are Jordan
disks.  Form the embedded dual $G$, with a vertex for each primal
face, and denote its exterior vertex by $0$.  Every face of $G$ is
a triangle with three distinct vertices, corresponding to a triple junction of the primal network. Each nonzero dual vertex inherits the label of its
corresponding primal face; distinct vertices may have the same label. Delete $0$
to obtain the finite-face dual $H$.

The graph $H$ is connected.  In the cyclic link of $0$, successive
nonzero neighbors are joined by the third edge of a triangular face;
these edges replace any path passing through $0$.  Every vertex of
$H$ was adjacent to $0$, so deleting $0$ exposes every vertex to
the same outer face.  Thus this embedding of $H$ is outerplanar.

There are no loops in $H$.  Two parallel finite--finite edges would
form a Jordan two-cycle.  Its side not containing $0$ cannot contain
a vertex: such a vertex would have to be adjacent to $0$, crossing
the cycle.  If it contains no vertex, that side cannot be filled by
triangular faces having three distinct vertices.  Thus $H$ is simple.
Every bounded face of $H$ remains a triangular face of $G$.

We relate the two degrees, retaining the possible parallel edges
from $0$.  For a finite vertex $x$ let $m_0(x)$ count its edges to
$0$, with multiplicity.  Then
\begin{equation}\label{eq:dual-degree}
 \deg_G(x)=\deg_H(x)+m_0(x).
\end{equation}
The link walk of $0$ becomes the outer boundary walk of $H$; its
occurrences of $x$ count exactly $m_0(x)$.  If $x$ occurred twice
and $H-x$ were connected, two distinct outer-face sectors at $x$
would separate nonempty groups of incident edges.  A simple path in
$H-x$ between neighbors in those two groups, together with their
edges to $x$, would form a Jordan cycle trapping one of these
outer-face sectors on its bounded side.  That sector could not belong
to the unbounded face.  Therefore
\begin{equation}\label{eq:one-exterior-side}
 x\text{ is not a cut vertex of }H\quad\Longrightarrow\quad m_0(x)=1.
\end{equation}
In particular a degree-one leaf represents a boundary
digon, and a non-cut degree-two vertex represents a boundary
triangle.  Initially every vertex of $H$ has degree at least two:
a degree-one vertex is not a cut vertex, so
\eqref{eq:dual-degree}--\eqref{eq:one-exterior-side} would give a
forbidden digon.  Since $H$ is finite, connected and has
minimum degree two, it contains a cycle and hence a bounded face.
All bounded faces are triangles.  Fix one such triangle as the core.
Its three vertices are pairwise adjacent and consequently have the
three distinct finite labels.

We describe a deletion procedure preserving these core vertices.
A block is a maximal two-connected subgraph or a bridge. The
block-cut tree has the blocks and cut vertices as its two classes of
vertices, with incidence as its edges; root it at the block
containing the core. A non-root leaf block has one attachment vertex,
namely its unique cut vertex shared with the rest of $H$.
An ear vertex of a triangulated polygon is a vertex of degree two.
Use this rooted tree, keeping the block containing the core.
If a non-core leaf block is a bridge, delete its far degree-one
endpoint.  A two-connected non-core leaf block is a triangulated
polygon with at most one attachment vertex.  For a polygon with
$N$ vertices, denote the numbers of triangles and interior diagonals
by $t,d$.  Euler's formula and triangle--edge incidence give
\[
 N-(N+d)+(t+1)=2,\qquad 3t=N+2d.
\]
The first equality gives $t-d=1$; substituting $d=t-1$ into the
second gives $3t=N+2t-2$.  Hence $t=N-2$ and $d=N-3$.
Its triangle-adjacency graph has one vertex for each triangle
and one edge for each shared interior diagonal. This graph is connected
(join interior points of two triangles by a path transverse to the
diagonals), with one fewer edge than vertices, and is therefore a
tree.  When $N>3$ it has at least two leaves, giving two distinct
degree-two ear vertices.  At least one differs from the attachment
vertex.  When $N=3$, either of the two non-attachment vertices can
be removed.  A chosen vertex has degree two in all of $H$ and is
not a cut vertex.  Delete it and recompute the blocks.  Degenerate
leaf blocks which become bridges are treated by the first rule.

After all non-core blocks have disappeared, the remaining block
is a triangulated polygon containing the core face.  Root its
triangle-adjacency tree at that face and remove the tip of any
non-root leaf triangle, where the tip is the vertex opposite the
edge shared with its parent. Such a tip is not a core vertex.  Indeed,
the triangles incident to a fixed polygon vertex form a connected
path in the triangle-adjacency tree.  If the leaf tip belonged to
the core, the path from that leaf to the core would consist of
triangles incident to the tip, whereas its first parent edge is
opposite the tip.  This is impossible.  Repeating removes precisely
the non-core vertices and leaves the core triangle.

Each graph deletion is realized by the corresponding geometric
deletion.  The primal face corresponding to a degree-one leaf is adjacent to
the exterior face and one finite face and, by \eqref{eq:one-exterior-side},
has exactly two sides;
apply Lemma~\ref{lem:digon-deletion}.  The primal face corresponding to a chosen degree-two ear is adjacent
to the exterior face and two distinct finite faces. The vertices of
$H$ corresponding to these two finite neighbors are adjacent in the
triangulation, so the primal face has exactly three sides; apply Lemma~\ref{lem:triangle-deletion}.
At such a face the three neighboring labels are distinct.  All
replacement arcs lie strictly inside the old face, extend the old
incident supporting curves, and preserve the neighboring pressure
differences.

These replacements preserve the required topology at
every step.  The deleted Jordan disk is divided into two or three
positive-area sectors, each attached to one distinct neighboring
primal face across its old open side.  Attaching a disk across
one boundary arc to a Jordan disk again gives a Jordan disk.  No
two old faces are identified, and no old face is split.  The
exterior gains exactly its attached sector and remains connected;
no exterior hole or crossing is introduced.  A triangle deletion
creates one regular $120$-degree triple junction, and a digon
deletion smoothly rejoins its two incident arcs.  The connected
primal graph changes locally by
\[
 (V,E,F)\longmapsto(V-2,E-3,F-1).
\]
Consequently the finite-face dual changes exactly from $H$ to
$H-x$, preserving connectedness and all three positive core faces.
These facts provide the induction hypotheses for the next
geometric deletion.

Every step decreases the number of finite faces by one.  Hence at
most $|V(H)|-3$ steps occur.  At the endpoint $H$ is the core
triangle.  The full dual has four vertices and all faces triangular,
so Euler's formula gives six edges.  Its three finite edges and
the required three exterior adjacencies already exhaust these edges;
thus it is $K_4$, and its embedded primal is also $K_4$.

If $D_t$ is the face removed at step $t$, with label $i$, write
$D_{t,j}$ for its sector assigned to neighboring label $j$, and set
$D_{t,j}=\varnothing$ for non-neighbors.  The neighboring labels
are distinct and the sectors partition $D_t$ up to their boundaries.
Thus the individual finite areas change by
\[
 v_i^{t+1}=v_i^t-|D_t|,\qquad
 v_j^{t+1}=v_j^t+|D_{t,j}|\quad(j\ne i,\ j\ne0).
\]
The removed component is non-core, while a core component of its
label remains, so $v_i^{t+1}>0$.  Summing the displayed equations,
and using $|D_t|=|D_{t,0}|+\sum_{j\ne0,i}|D_{t,j}|$, gives exactly
\[
 A_{t+1}=A_t-|D_{t,0}|,\qquad |D_{t,0}|>0.
\]
All other sectors remain finite.  Thus total finite area strictly
decreases whenever a step is performed.
\end{proof}

\begin{proposition}[Excluding a comparison minimum below one]
\label{prop:no-strict-sublevel}
There is no comparison minimum under a perimeter bound with $0<Q<1$.
\end{proposition}
\begin{proof}
Choose an actual fixed-area minimizing cluster $C_*$ representing
such a comparison minimum, as in Proposition~\ref{prop:attainment}.
Lemma~\ref{lem:disk-faces} and Lemma~\ref{lem:no-digon} give a
compact connected regular network with Jordan-disk finite faces and
no digon.  Proposition~\ref{prop:components} gives connected
exterior and the component restrictions.  Proposition~\ref{prop:two-types}
leaves the 4-daisy or the all-boundary type.  The former admits
the strict, perimeter-nonincreasing improvement of Proposition~\ref{prop:no-daisy},
so every finite face meets the exterior.

Apply Lemma~\ref{lem:pruning}.  If no deletion occurs, $C_*$ is
already a $K_4$ network and hence standard by Lemma~\ref{lem:k4},
contradicting $Q(C_*)<1$.  Otherwise the endpoint is a compact
equilibrium $K_4$ with the same actual pressures $p$ and strictly
smaller total finite area.  Lemma~\ref{lem:k4} identifies it as
standard.  The pressure realization of the standard family then
implies $p_i>1$ for all $i$.  Together with
Lemma~\ref{lem:capped-pressure}, this gives
\[
 1<p_i<\kappa_i\quad(1\le i\le3),
\]
where $\kappa$ denotes the standard pressure vector at the areas of
$C_*$.  Let $A_{\mathrm{std}}$ be total finite area in pressure
coordinates.  The whole segment
$\kappa(t)=\kappa+t(p-\kappa)$ lies in $(1,\infty)^3$.
The standard-family identity $dA_{\mathrm{std}}=-\sum_iw_{i0}\,d\kappa_i$
therefore yields
\[
 \begin{split}
 A_{\mathrm{std}}(p)-A_{\mathrm{std}}(\kappa)
 &=\int_0^1\sum_{i=1}^3
        w_{i0}(\kappa(t))(\kappa_i-p_i)\,dt>0.
 \end{split}
\]
Every weight is positive.  But $A_{\mathrm{std}}(\kappa)=A(C_*)$
and the pruning endpoint has area $A_{\mathrm{std}}(p)<A(C_*)$.
This is the opposite inequality, a contradiction.
\end{proof}

\begin{proof}[Proof of Theorem~\ref{thm:main}]
If some finite-perimeter competitor at a positive area vector had
$P<T(v)$, its positive comparison ratio would be below one.
Fix a perimeter bound at least as large as its perimeter.
Proposition~\ref{prop:attainment} would produce an attained strict
sublevel minimum, contradicting Proposition~\ref{prop:no-strict-sublevel}.
Thus $P\ge T(v)$ for every positive area vector.  Since the standard
cluster is admissible, the true profile equals $T$.

Let now $P(C)=T(v)$.  The cluster is an actual fixed-area minimizer,
so the geometric inputs apply.  The global inequality just proved
gives $Q\ge1$ for all positive area vectors; hence $C$ is an
unconstrained global minimum of the comparison ratio.  The strict
relabeling arguments of Proposition~\ref{prop:components} and
Proposition~\ref{prop:no-daisy} apply also at $Q=1$.
Together with Proposition~\ref{prop:two-types} and the absence of
digons, they imply that every finite component meets the exterior.

Both signs of each sufficiently small compactly supported smooth
variation remain admissible, with the three areas positive.  At
$Q=1$ its derivative is
\[
 0=DQ=\frac{DP-\sum_i\kappa_iDv_i}{S(v)}.
\]
Stationarity gives $DP=\sum_ip_iDv_i$, and the area differential is
onto by Lemma~\ref{lem:area-surjectivity}.  It follows that $p=\kappa$.
If the network were not $K_4$, Lemma~\ref{lem:pruning} would give
a compact standard $K_4$ endpoint at these same pressures and with
strictly smaller total finite area.  Uniqueness of pressure
realization identifies this endpoint with the original standard
model at area vector $v$, whose total area is $A(v)$.  The strict
decrease is impossible.  Thus the original network is $K_4$ and is
standard by Lemma~\ref{lem:k4}.  The volume-coordinate theorem for
the standard family gives uniqueness at the prescribed labeled
areas, modulo a hyperbolic isometry and null-set changes.
\end{proof}

\begin{corollary}[Curvature rescaling]\label{cor:curvature-rescaling}
For curvature $-a^2$, $a>0$, the same uniqueness assertion holds and
the standard perimeter profile is
\[
 T_{-a^2}(v)=a^{-1}T_{-1}(a^2v).
\]
\end{corollary}
\begin{proof}
Multiplying the curvature-$-1$ metric by $a^{-2}$ multiplies areas
by $a^{-2}$ and lengths by $a^{-1}$.  This constant rescaling gives
a bijection of competitors and equality cases, proving both claims.
\end{proof}

\section*{Acknowledgment of assistance}
The author used ChatGPT in preparing and revising this manuscript and takes responsibility for its mathematical content and references.

\end{document}